\documentclass[11pt]{article}
\usepackage[margin=1in]{geometry}
\usepackage[T1]{fontenc}
\usepackage{amsmath,amssymb,amsthm,mathtools,enumitem,booktabs}
\usepackage[dvipsnames]{xcolor}
\usepackage[colorlinks=true,linkcolor=blue,citecolor=blue,urlcolor=blue]{hyperref}

\newtheorem{theorem}{Theorem}[section]
\newtheorem{lemma}[theorem]{Lemma}
\newtheorem{proposition}[theorem]{Proposition}
\newtheorem{corollary}[theorem]{Corollary}
\theoremstyle{remark}
\newtheorem{remark}[theorem]{Remark}

\newcommand{\R}{\mathbb{R}}
\newcommand{\Sph}{\mathbb{S}}

\newcommand{\Hh}{\mathcal{H}}
\newcommand{\X}{\mathcal{X}}
\newcommand{\Rp}{\mathbb{R}^N_+}
\newcommand{\dd}{\,\mathrm{d}}

\title{Positive self-similar solutions of a parabolic equation with\\
spatially weighted absorption in a half-space}
\author{Waad Al Sayed\thanks{Department of Mathematics and Statistics, College of Sciences and Humanities, Fahad Bin Sultan University, P.O. Box 15700, Tabuk 71454, Saudi Arabia. \texttt{wsaid@fbsu.edu.sa}}
\and
Mustapha Jazar\thanks{Lebanese University, LaMA-Liban, Tripoli, Lebanon. \texttt{mjazar@laser-lb.org}}}
\date{}

\begin{document}
\maketitle

\begin{abstract}
We study positive solutions of $u_t-\Delta u+x_1^{2\alpha}|u|^{q-1}u=0$ in the
half-space $\{x_1>0\}$, with $\alpha>-1$ and $q>1$, vanishing on $\{x_1=0\}$ and
at $t=0$ away from the origin. Let $q_c=(N+3+2\alpha)/(N+1)$. For $1<q<q_c$ we prove
that there is exactly one positive self-similar solution continuous away from the origin,
and that its profile satisfies Gaussian bounds that are sharp up to multiplicative constants.
For $q\ge q_c$ there is none. We also show that, after removing the Gaussian factor, the profile has an angular limit shape $\Phi$ in every mode, give an integral formula for its Kummer coefficients, prove that $\Phi$ is real-analytic and is the uniform limit of the rescaled profiles, that $\Phi/(\eta_1/|\eta|)$ is monotone along meridians, and that $\Phi$ is proportional to $\eta_1/|\eta|$ if and only if $2\alpha+q=1$. Every self-similar solution of finite energy, or continuous away from the origin, is dominated in absolute value by the positive one, and for $(N+4+2\alpha)/(N+2)\le q<q_c$ the only nontrivial ones are $\pm u_s$. For $\alpha=0$ existence, uniqueness and the upper bound are due to Marcus and V\'eron; the lower bound with the sharp exponent and the description of $\Phi$ are new also in that case. The existence and uniqueness proofs use a Brezis--Oswald type comparison principle and explicit barriers, given by explicit formulas; the analysis of $\Phi$ uses Kummer functions.
\medskip

\end{abstract}

%---------------------------------------------------------------------
\section{Introduction}
Let $N\ge 2$, $\alpha>-1$, $q>1$, $\Rp=\{x\in\R^N:\ x_1>0\}$ and $E=\Rp\times(0,\infty)$.
We consider positive solutions of
\begin{equation}\label{eq:pde}
u_t-\Delta u+x_1^{2\alpha}|u|^{q-1}u=0\quad\text{in }E,
\end{equation}
\begin{equation}\label{eq:bc}
u(x,t)=0\ \text{ if } x_1=0,\ t>0,\qquad u(x,t)\to0\ \text{ as } t\to0,\ x\neq0 .
\end{equation}
Equation \eqref{eq:pde} is invariant under $u\mapsto u_\mu(x,t)=\mu^{2\lambda}u(\mu x,\mu^2t)$,
where
\[
\lambda=\frac{1+\alpha}{q-1},\qquad q_c=\frac{N+3+2\alpha}{N+1}.
\]
We look for self-similar solutions $u(x,t)=t^{-\lambda}h(x/\sqrt t)$. Put $K(\eta)=e^{|\eta|^2/4}$, $r=|\eta|$, $\sigma=\eta/r$, $\sigma_1=\eta_1/r$,
$s=2\lambda-N$. Then $h$ solves
\begin{equation}\label{eq:E}
L_Kh-\lambda h+\eta_1^{2\alpha}h^q=0\ \text{ in }\Rp,\qquad h=0\ \text{ on }\partial\Rp,
\qquad L_Kv:=-K^{-1}\operatorname{div}(K\nabla v)=-\Delta v-\tfrac12\,\eta\cdot\nabla v .
\end{equation}

\begin{theorem}\label{thm:main}
Let $N\ge2$, $\alpha>-1$, $q>1$. Let $\Hh$ be the energy class of Section~\ref{sec:prelim} and
$\mathcal S$ the set of continuous positive profiles of Section~\ref{sec:cont}
(continuous on $\overline{\Rp}$, vanishing on $\partial\Rp$, with $|\eta|^{2\lambda}h\to0$).
\begin{enumerate}[label=\textup{(\roman*)},leftmargin=*]
\item If $q\ge q_c$, then $\mathcal S=\emptyset$ and $\Hh=\{0\}$.
\item If $1<q<q_c$, then $\mathcal S=\Hh\setminus\{0\}=\{h\}$ for a single function $h$. It
belongs to $C^2(\Rp)\cap C(\overline{\Rp})$, is positive in $\Rp$, and there are
$c,C,R_1>0$ with
\begin{equation}\label{eq:sharp}
c\,\eta_1|\eta|^{2\lambda-N-1}e^{-|\eta|^2/4}\le h(\eta)\le C\,\eta_1|\eta|^{2\lambda-N-1}e^{-|\eta|^2/4}
\qquad(|\eta|\ge R_1).
\end{equation}
\item For $q<q_c$, $u_s(x,t)=t^{-\lambda}h(x/\sqrt t)$ is a positive solution of
\eqref{eq:pde}--\eqref{eq:bc}, continuous on $\overline E\setminus\{(0,0)\}$, and it is the
\emph{only} positive self-similar solution in $C(\overline E\setminus\{(0,0)\})$. Moreover
$\int_{\Rp} x_1u_s(x,t)\dd x=t^{\frac{N+1}2-\lambda}\int_{\Rp}\eta_1h\dd\eta\to\infty$ as $t\to0$.
For $q\ge q_c$ there is no such solution.
\end{enumerate}
\end{theorem}

For $\alpha=0$, parts (i)--(iii) of Theorem~\ref{thm:main}, except the lower bound in \eqref{eq:sharp} and the statement on the mass of $x_1u_s$, are contained in \cite[Theorem~2.1 and Lemma~2.5]{MV01}; see the comments below.

Solutions of \eqref{eq:pde} need not be positive, and the equation is odd. Let $\Sigma$ be the set of
finite-energy self-similar profiles of any sign, i.e.\ $w\in\X$ satisfying \eqref{eq:weak} with $h^q$ replaced by
$|w|^{q-1}w$ (Section~\ref{sec:signed}); $\Hh=\Sigma\cap\{w\ge0\}$. Put
\[
q_2=\frac{N+4+2\alpha}{N+2}\in(1,q_c),\qquad\text{so that}\quad \lambda\le\lambda_2:=\frac{N+2}2\iff q\ge q_2 .
\]

\begin{theorem}\label{thm:signed}
Let $N\ge2$, $\alpha>-1$, $q>1$.
\begin{enumerate}[label=\textup{(\roman*)},leftmargin=*]
\item The profile of every self-similar solution of \eqref{eq:pde}--\eqref{eq:bc} in $C(\overline E\setminus\{(0,0)\})$, of any sign, belongs to $\Sigma$. If $q\ge q_c$, then $\Sigma=\{0\}$.
\item If $1<q<q_c$, then $|w|\le h$ in $\Rp$ for every $w\in\Sigma$; hence $|u|\le u_s$ for every self-similar solution $u\in C(\overline E\setminus\{(0,0)\})$.
\item If $w\in\Sigma$ changes sign, then $q<q_2$. Consequently, if $q_2\le q<q_c$, then $\Sigma=\{0,h,-h\}$, and $u_s$, $-u_s$ are the only nontrivial self-similar solutions in $C(\overline E\setminus\{(0,0)\})$.
\end{enumerate}
\end{theorem}

\paragraph{Relation with the work of Marcus and V\'eron ($\alpha=0$).}
For $\alpha=0$ the half-space problem was solved by Marcus and V\'eron \cite[Theorem~2.1]{MV01} (announced as Theorem~0.4 there; see also the notes \cite{MV96,MV97}).
For $1<q<q_c$ they prove that there is exactly one positive self-similar solution $U_s\in C(\overline E\setminus\{(0,0)\})$, that it dominates every solution in this class \cite[Theorem~2.1(iii) and Lemma~2.7]{MV01}, that its profile satisfies the upper bound in \eqref{eq:sharp}, and that it satisfies the lower bound $c_m\eta_1|\eta|^{2m-N-1}e^{-|\eta|^2/4}$ for every $m<\lambda$, with $c_m$ depending on $m$. They also prove that for $(N+4)/(N+2)<q<q_c$ the only self-similar solutions are $0$ and $\pm U_s$ \cite[Lemma~2.6]{MV01}, and that isolated singularities on the lateral boundary are removable for $q\ge q_c$ \cite[Proposition~0.1]{MV01}. Their uniqueness class is the same as here: a self-similar solution continuous on $\overline E\setminus\{(0,0)\}$ has a profile of finite energy \cite[Lemma~2.5]{MV01}. The Dirichlet spectrum of $L_K$ on the half-space is computed in \cite[Lemma~2.2]{MV01}, following \cite{EK}. Positive solutions with an isolated singularity on the lateral boundary are classified, and shown to converge to $U_s$ after rescaling, in \cite[Theorem~0.5]{MV01}; signed solutions are treated in \cite{MV01b}.

Building on the foundational results of Marcus and V\'eron \cite{MV01} for the unweighted half-space problem, the present study investigates self-similar solutions in the presence of the boundary-normal absorption weight $x_1^{2\alpha}$, with $\alpha>-1$. This setting includes singular weights when $-1<\alpha<0$ and requires additional care in the boundary estimates and weighted variational arguments. The analysis also yields refined information on the asymptotic behaviour of the self-similar profile, including matching algebraic exponents in the two-sided Gaussian bounds. Furthermore, the far-field angular profile is characterized through Kummer-function representations, together with results on its analyticity, uniform angular convergence, and monotonicity properties. These findings complement the existing theory and provide a more detailed description of the qualitative and asymptotic structure of self-similar solutions.

\paragraph{What is new here.}
(a) The weight $x_1^{2\alpha}$ with $\alpha>-1$, including singular weights $\alpha<0$; these change the boundary analysis, which is why the test-function and Hardy arguments of Lemma~\ref{lem:test}(c) and Section~\ref{sec:cont} are written out.
(b) The lower bound in \eqref{eq:sharp} has the same exponent as the upper bound; this is new also for $\alpha=0$.
(c) The proofs use explicit barriers, among them the exact solution $W=M\eta_1^{-2\lambda}$ (Lemma~\ref{lem:W}), and no information on linear ODE solutions. The function $z_1\sigma_1$ is a supersolution only outside a large ball; Marcus and V\'eron use it there (Remark~\ref{rem:nosuper}).
(d) Theorem~\ref{thm:signed}(iii) includes the endpoint $q=q_2$ and allows singular weights; its proof uses only $w^\pm$ and the second eigenvalue.
(e) Section~\ref{sec:asym} is new also for $\alpha=0$: the angular limit $\Phi$ exists mode by mode (by an explicit reduction to Kummer's equation), it depends only on $\sigma_1$, its coefficients are given by integral formulas, it is proportional to $\sigma_1$ exactly when $2\alpha+q=1$, and, for every $q<q_c$, $\Phi/\sigma_1$ is monotone with $\operatorname{sign}c_3=-\operatorname{sign}(2\alpha+q-1)$; moreover $H(r,\cdot)\to\Phi$ uniformly, $\Phi$ is real-analytic, and near $q_c$ the sign of every coefficient follows a $\Gamma$-function law, except at the exceptional values discussed after Theorem~\ref{thm:nearcrit}.
The uniqueness among continuous profiles (Section~\ref{sec:cont}) and the domination $|w|\le h$ (Theorem~\ref{thm:signed}(ii)) follow the strategy of \cite{MV01}; we give complete proofs because of the weight. Domination of \emph{all} solutions, not only self-similar ones, by $u_s$ is proved in \cite[Lemma~2.7]{MV01} for $\alpha=0$ and is not addressed here.

\paragraph{Other literature.}
For the unweighted equation in $\R^N$, measure initial data and isolated initial singularities go back to \cite{BF,BP}; the very singular self-similar solution was constructed in \cite{BPT}, obtained as the limit of source-type solutions in \cite{KP}, and isolated positive singularities were classified in \cite{Os}. Initial traces for absorption equations were developed in \cite{MV99,MV02}, in analogy with the boundary trace of elliptic equations \cite{MV03}; the absorption weighted \emph{in time}, $\partial_tu-\Delta u+t^\alpha|u|^{q-1}u=0$, is treated in \cite{MV02} for $\alpha\ge0$ and in \cite{ASV} for $\alpha>-1$, in $\R^N$ and in domains with compact boundary. There the very singular solution has the same exponent $\lambda=(1+\alpha)/(q-1)$ as here, but the critical exponent is $1+2(1+\alpha)/N$, because the threshold is $\lambda>N/2$ instead of $\lambda>\lambda_1=(N+1)/2$. Hisa, Ishige and Takahashi \cite{HIT} study the source equation $u_t=\Delta u+u^p$ in a half-space with a Dirichlet condition; their initial trace is that of $x_Nu$, the quantity whose mass blows up in Theorem~\ref{thm:main}(iii). Self-similar solutions with an isolated singularity exist here because a half-space is scale invariant (cones are treated for $\alpha=0$ in \cite[Theorem~2.8]{MV01}); for domains with compact boundary there are none. We do not claim any statement on initial traces of general solutions.

%---------------------------------------------------------------------
\section{Preliminaries}\label{sec:prelim}

\subsection{Reduction}
For $u(x,t)=t^{-\lambda}h(\eta)$, $\eta=x/\sqrt t$, one has
$u_t=t^{-\lambda-1}(-\lambda h-\frac12\eta\cdot\nabla h)$, $\Delta u=t^{-\lambda-1}\Delta h$ and,
since $\lambda q=\lambda+1+\alpha$,
$x_1^{2\alpha}u^q=t^{-\lambda-1}\eta_1^{2\alpha}h^q$. Hence \eqref{eq:pde} is equivalent to
\eqref{eq:E}. For $v\ge0$ we write $\mathcal M(v)=L_Kv-\lambda v+\eta_1^{2\alpha}v^q$.

\subsection{Function spaces}
$L^2_K(\Rp)$ is $L^2(\Rp;K\dd\eta)$; $W^{1,2}_{0,K}(\Rp)$ is the closure of $C_c^\infty(\Rp)$ for
$\|v\|^2=\int_{\Rp}(|\nabla v|^2+v^2)K\dd\eta$; $L^{q+1}_{\eta_1^{2\alpha}K}(\Rp)=L^{q+1}(\Rp;\eta_1^{2\alpha}K\dd\eta)$. Put
\[
\X=W^{1,2}_{0,K}(\Rp)\cap L^{q+1}_{\eta_1^{2\alpha}K}(\Rp).
\]
By H\"older, $\int\eta_1^{2\alpha}|h|^{q}|\varphi|K<\infty$ for $h,\varphi\in\X$. Let $\Hh$ be the set of $h\in\X$ with $h\ge0$ satisfying
\begin{equation}\label{eq:weak}
\int_{\Rp}\nabla h\cdot\nabla\varphi\,K\dd\eta+\int_{\Rp}\eta_1^{2\alpha}h^q\varphi\,K\dd\eta
=\lambda\int_Ah\varphi\,K\dd\eta\qquad\text{for all }\varphi\in\X .
\end{equation}

\begin{lemma}\label{lem:reg}
Every $h\in\Hh$ is bounded, and $\sup_{B_1(y)\cap\Rp}h\le C\,e^{-(|y|-1)_+^2/8}\,\|h\|_{L^2_K(\Rp)}$ with $C=C(N,\lambda)$.
In particular $h(\eta)\to0$ as $|\eta|\to\infty$.
\end{lemma}
\begin{proof}
For $w\in C_c^\infty(\Rp)$ one has $\int|\nabla(K^{-1/2}w)|^2K=\int|\nabla w|^2+(\frac{|\eta|^2}{16}+\frac N4)w^2$
(integrate $-\frac12 w\,\eta\cdot\nabla w$ by parts). Hence $v=K^{1/2}h\in W^{1,2}_0(\Rp)$ and, by polarisation,
\eqref{eq:weak} reads
\[
\int\nabla v\cdot\nabla\tilde v+\Big(\tfrac{|\eta|^2}{16}+\tfrac N4-\lambda\Big)v\tilde v+\int\eta_1^{2\alpha}K^{\frac{1-q}2}v^q\tilde v=0,\qquad\tilde v=K^{1/2}\varphi,\ \varphi\in\X .
\]
Take $\tilde v=\zeta\,T_\tau(v)$ with $0\le\zeta\in C_c^\infty(\R^N)$ and $T_\tau(t)=\min(t/\tau,1)$; this is admissible because
$\tilde v\le\|\zeta\|_\infty v/\tau$. Dropping the nonnegative terms and letting $\tau\to0$ shows that the zero extension $\hat v$
of $v$ satisfies $-\Delta\hat v\le\lambda\hat v$ in $\mathcal D'(\R^N)$. Moser iteration for this constant-coefficient inequality \cite[Theorem~8.17]{GT} gives
$\sup_{B_1(y)}\hat v\le C(N,\lambda)\|\hat v\|_{L^2(B_2(y))}$, with $C$ independent of $y$. Finally $h=K^{-1/2}v\le e^{-(|y|-1)_+^2/8}\sup_{B_1(y)}v$
on $B_1(y)$, and $\|\hat v\|_{L^2(B_2(y))}\to0$ as $|y|\to\infty$ because $\hat v\in L^2(\R^N)$.
\end{proof}

\subsection{The bottom of the spectrum}
Let $\psi_1(\eta)=\eta_1e^{-|\eta|^2/4}$. Then $L_K\psi_1=\lambda_1\psi_1$ with $\lambda_1=\frac{N+1}2$,
and $\psi_1\in\X$, since $\int\eta_1^{2\alpha}\psi_1^{q+1}K=\int\eta_1^{2\alpha+q+1}e^{-q|\eta|^2/4}<\infty$
(here $2\alpha+q+1>0$ is used). Note that
\begin{equation}\label{eq:qc}
\lambda>\lambda_1\iff q<q_c,\qquad s=2\lambda-N>1\iff q<q_c .
\end{equation}

\begin{lemma}\label{lem:spec}
\textup{(i)} $\int_{\Rp}|\nabla\varphi|^2K\ge\lambda_1\int_{\Rp}\varphi^2K$ for all $\varphi\in W^{1,2}_{0,K}(\Rp)$, and the
infimum of the Rayleigh quotient over $C_c^\infty(\Rp)$ is $\lambda_1$. Consequently, if
$D_n=\Rp\cap B_n(0)$ and $\lambda_1(D_n)$ is the first Dirichlet eigenvalue of $L_K$ in $D_n$,
then $\lambda_1(D_n)\downarrow\lambda_1$.

\textup{(ii)} If $\varphi\in W^{1,2}_{0,K}(\Rp)$ and $\int_{\Rp}\varphi\psi_1K=0$, then $\int_{\Rp}|\nabla\varphi|^2K\ge\lambda_2\int_{\Rp}\varphi^2K$, $\lambda_2=\frac{N+2}2$.
\end{lemma}
\begin{proof}
Write $\varphi=\psi_1v$ with $v\in C_c^\infty(\Rp)$. Since
$|\nabla\varphi|^2=\nabla\psi_1\cdot\nabla(\psi_1v^2)+\psi_1^2|\nabla v|^2$, integration by parts gives
$\int|\nabla\varphi|^2K-\lambda_1\int\varphi^2K=\int\psi_1^2|\nabla v|^2K\ge0$. Cutting off $\psi_1$
gives test functions with quotient arbitrarily close to $\lambda_1$. Monotonicity and exhaustion give the last claim; the inequality extends to $W^{1,2}_{0,K}(\Rp)$ by density.

(ii) This is \cite[Lemma~2.2]{MV01} (following \cite{EK}); we recall the argument. For $\varphi\in C_c^\infty(\Rp)$ and $w=K^{1/2}\varphi$,
$\int|\nabla\varphi|^2K=\int|\nabla w|^2+(\frac{|\eta|^2}{16}+\frac N4)w^2$ and $\int\varphi^2K=\int w^2$ (proof of Lemma~\ref{lem:reg}). So the quadratic form of $L_K$ on $L^2_K(\Rp)$ is unitarily equivalent to that of $\mathcal O+\frac N4$ on $L^2(\Rp)$, where $\mathcal O=\sum_j(-\partial_j^2+\frac{\eta_j^2}{16})$ with the Dirichlet condition on $\{\eta_1=0\}$, and the identity extends to the form domains by closure. The functions $H_k(t/2)e^{-t^2/8}$ ($H_k$ Hermite) are a complete set of eigenfunctions of $-\frac{d^2}{dt^2}+\frac{t^2}{16}$ in $L^2(\R)$, with eigenvalues $\frac{2k+1}4$; those with $k$ odd vanish at $0$ and are complete in $L^2(0,\infty)$. Hence products of these functions, with $k_1$ odd, form a complete eigenbasis of $\mathcal O$, and the spectrum of $L_K$ is $\{\frac N2+\frac{k_1+\dots+k_N}2:\ k_1\text{ odd}\}=\{\frac{N+k}2\}_{k\ge1}$. The lowest level ($k_1=1$, other $k_j=0$) is simple, with eigenfunction $K^{-1/2}\eta_1e^{-|\eta|^2/8}=\psi_1$. The min--max principle gives (ii).
\end{proof}

\subsection{A comparison principle}
For $Q\subset\Rp$ open and $u\ge0$ locally bounded, set
$\mathcal B(u;\varphi)=\int_Q\nabla u\cdot\nabla\varphi\,K+\int_Q\eta_1^{2\alpha}u^q\varphi K-\lambda\int_Qu\varphi K$.

\begin{lemma}[comparison]\label{lem:comp}
Let $Q\subset\Rp$ be open and let $u_1,u_2\in W^{1,2}_{\rm loc}(Q)$, $u_1\ge0$ bounded,
$u_2>0$ a.e. Assume:
\begin{enumerate}[label=\textup{(H\arabic*)},leftmargin=*]
\item $\mathcal B(u_1;\varphi)\le0\le\mathcal B(u_2;\varphi)$ for every $0\le\varphi\in W^{1,2}_0(Q)\cap L^\infty$
with bounded support, vanishing outside $\{u_1>u_2+\varepsilon\}$ for some $\varepsilon>0$, and
such that the integrals converge;
\item $u_1\in L^2_K(Q)$, and $\eta_1^{2\alpha}u_1^{q+1}K\in L^1(S_\varepsilon)$ for every $\varepsilon>0$ (automatic if $S_\varepsilon$ is relatively compact in $\Rp$), where $S_\varepsilon$ is defined in \textup{(H3)};
\item for every $\varepsilon>0$ the set $S_\varepsilon=\{u_1>u_2+\varepsilon\}$ is bounded,
$\nabla u_1,\nabla u_2\in L^2(S_\varepsilon)$, and $(u_1-u_2-\varepsilon)_+\in W^{1,2}_0(Q)$.
\end{enumerate}
Then $u_1\le u_2$ a.e.\ in $Q$.
\end{lemma}
\begin{proof}
This is the argument of Brezis--Oswald \cite{BO} and D\'iaz--Saa \cite{DS}; \cite[Lemma~2.3]{MV01} is a version for $\alpha=0$ in exterior domains. Fix $\varepsilon\in(0,1)$, put $a=u_1+\varepsilon/2$,
$b=u_2+\varepsilon$, and
\[
w_1=\frac{(a^2-b^2)_+}{a},\qquad w_2=\frac{(a^2-b^2)_+}{b}.
\]
Both are nonnegative, bounded and supported in $S_{\varepsilon/2}$. Writing $f=(u_1-u_2-\frac\varepsilon2)_+$, we have $w_1=f(1+\frac ba)$, $w_2=f(1+\frac ab)$, where the second factors are bounded on $S_{\varepsilon/2}$ with gradients in $L^2(S_{\varepsilon/2})$ by (H3); so $w_i\in W^{1,2}(Q)$, $0\le w_i\le C_\varepsilon f$, and $w_i\in W^{1,2}_0(Q)$ because $f\in W^{1,2}_0(Q)$ (if $0\le u\le v$ with $u\in W^{1,2}$, $v\in W^{1,2}_0$, then $u\in W^{1,2}_0$).
On the support $a<2u_1$, so $w_1\le 2u_1$ and $u_2^qw_2\le4u_1^{q+1}$; by (H2) all integrals converge.
Testing with $w_1$ and $w_2$ and subtracting,
\[
\int(\nabla u_1\!\cdot\!\nabla w_1-\nabla u_2\!\cdot\!\nabla w_2)K
+\int\eta_1^{2\alpha}(u_1^qw_1-u_2^qw_2)K\le\lambda\int(u_1w_1-u_2w_2)K .
\]
A direct computation on $S_{\varepsilon/2}$ gives the identity
\[
\nabla u_1\!\cdot\!\nabla w_1-\nabla u_2\!\cdot\!\nabla w_2
=\Big|\nabla u_1-\tfrac ab\nabla u_2\Big|^2+\Big|\nabla u_2-\tfrac ba\nabla u_1\Big|^2\ge0 .
\]
Moreover $u_1w_1-u_2w_2=\big(\tfrac{u_1}a-\tfrac{u_2}b\big)(a^2-b^2)_+$ is bounded in modulus by
$a^2\mathbf 1_{S_{\varepsilon/2}}\le4u_1^2\in L^1(K)$, and tends to $0$ a.e.\ as $\varepsilon\to0$
(this uses $u_2>0$): at a point with $u_1>u_2>0$ the bracket tends to $0$, and at a point with
$u_1\le u_2$ the integrand vanishes. Hence the right-hand side tends to $0$. On the other hand
\[
u_1^qw_1-u_2^qw_2=\Big(\tfrac{u_1^q}{a}-\tfrac{u_2^q}{b}\Big)(a^2-b^2)_+\ge0,
\]
because on $S_{\varepsilon/2}$ we have $u_1>u_2$ and $t\mapsto t^q/(t+c)$ is increasing, with
$u_2+\varepsilon/2<u_2+\varepsilon$. By Fatou's lemma,
\[
\int_{\{u_1>u_2\}}\eta_1^{2\alpha}\big(u_1^{q-1}-u_2^{q-1}\big)\big(u_1^2-u_2^2\big)K\le0,
\]
and the integrand is positive on $\{u_1>u_2\}$. Hence $|\{u_1>u_2\}|=0$.
\end{proof}

\begin{remark}\label{rem:H1}
Lemma~\ref{lem:test} below shows how (H1) is verified in the applications: for classical sub- or supersolutions $u$ whose sets $S_\varepsilon$ are relatively compact in $\Rp$, by part~(a);
for positive classical supersolutions with $(L_K-\lambda)u\ge0$ (the barriers $\Phi$, $C_2U$, $V_\varepsilon$), by part~(b), even when the support of the test function reaches $\partial\Rp$;
and for weak solutions in the sense of \eqref{eq:weak} ($h_n$, $h$, $h'$) directly, because the test functions $w_i$ are bounded by a multiple of $u_1$ and therefore belong to $\X$ (resp.\ $X_n$).
\end{remark}

\begin{lemma}[test functions]\label{lem:test}
Write $\mathcal B(u;\varphi)=\int_{\Rp}\nabla u\cdot\nabla\varphi\,K+\int_{\Rp}\eta_1^{2\alpha}u^q\varphi K-\lambda\int_Au\varphi K$.
\begin{enumerate}[label=\textup{(\alph*)},leftmargin=*]
\item Let $Q\subset\Rp$ be open (e.g.\ $Q=\Rp$) and $h\in C^2(Q)$, $h\ge0$. For every $\varphi\in W^{1,2}(Q)$ with compact support in $Q$,
\[
\mathcal B(h;\varphi)=\int_{\Rp}\mathcal M(h)\,\varphi\,K .
\]
In particular $\mathcal B(h;\varphi)\le0$ (resp.\ $\ge0$, $=0$) for all such $\varphi\ge0$ if $\mathcal M(h)\le0$ (resp.\ $\ge0$, $=0$) pointwise.
\item Let $Q\subset\Rp$ be open and $u\in C^2(Q)$, $u>0$, such that $u$, $\nabla u$ and $(L_K-\lambda)u$ are bounded on $Q\cap B_\rho$ for every $\rho$, and $(L_K-\lambda)u\ge0$ in $Q$.
Then $\mathcal B(u;\varphi)\ge0$ for every $0\le\varphi\in W_0^{1,2}(Q)$ with bounded support (the value can be $+\infty$ only through the nonnegative term $\int\eta_1^{2\alpha}u^q\varphi K$).
\item Let $h\in C^2(\Rp)\cap\X$, $h\ge0$, with $\mathcal M(h)=0$. Then \eqref{eq:weak} holds for all $\varphi\in\X$.
\end{enumerate}
\end{lemma}
\begin{proof}
(a) Extend $\varphi$ by $0$; then $\varphi\in W^{1,2}(\R^N)$ with support in a compact $F\subset\Rp$. Let $d=\operatorname{dist}(F,\partial\Rp)>0$ and $F'=\{\operatorname{dist}(\cdot,F)\le d/2\}$, a compact subset of $\Rp$.
For $\varepsilon<d/2$ the mollifications $\varphi_\varepsilon\in C_c^\infty(\Rp)$ are supported in $F'$, satisfy $\varphi_\varepsilon\ge0$ if $\varphi\ge0$, and $\varphi_\varepsilon\to\varphi$ in $W^{1,2}(\R^N)$.
Since $h\in C^2(\Rp)$ and $L_Kh=-K^{-1}\operatorname{div}(K\nabla h)$, integration by parts gives
$\int\nabla h\cdot\nabla\varphi_\varepsilon K=\int(L_Kh)\varphi_\varepsilon K$, that is, $\mathcal B(h;\varphi_\varepsilon)=\int\mathcal M(h)\varphi_\varepsilon K$.
On $F'$ the functions $K\nabla h$, $K\eta_1^{2\alpha}h^q$, $Kh$ and $K\mathcal M(h)$ are bounded (continuous on a compact set), so every term converges as $\varepsilon\to0$, using $\varphi_\varepsilon\to\varphi$ in $L^2(F')\subset L^1(F')$ and $\nabla\varphi_\varepsilon\to\nabla\varphi$ in $L^2(F')$.

(b) Let $\varphi\in W_0^{1,2}(Q)$, $\varphi\ge0$, with support in $B_\rho$. Choose $\varphi_i\in C_c^\infty(Q)$ with $\varphi_i\to\varphi$ in $W^{1,2}$ and replace $\varphi_i$ by $\varphi_i^+\zeta$, where $\zeta$ is a fixed cut-off equal to $1$ on $B_\rho$ and supported in $B_{\rho+1}$. These functions are Lipschitz, nonnegative, compactly supported in $Q$, supported in $B_{\rho+1}$,
are bounded in $W^{1,2}$ and converge to $\varphi$ in $L^2$, hence weakly in $W^{1,2}$. As in (a), $\int\nabla u\cdot\nabla\varphi_i^+\zeta K-\lambda\int u\varphi_i^+\zeta K=\int(L_K-\lambda)u\,\varphi_i^+\zeta K\ge0$. The hypotheses on $u$ allow us to pass to the limit in the left-hand side; hence the linear part of $\mathcal B(u;\varphi)$ is $\ge0$, and
$\mathcal B(u;\varphi)=\text{(linear part)}+\int\eta_1^{2\alpha}u^q\varphi K\ge0$.

(c) Let $\chi_k(r)$ be a cut-off ($\chi_k=1$ for $r\le k$, $0$ for $r\ge k+1$, $|\nabla\chi_k|\le2$) and $\theta_j(t)=\min\{1,(jt-1)_+\}$, so $\theta_j\uparrow1$ on $(0,\infty)$ and $0\le\theta_j'\le j\mathbf 1_{(1/j,2/j)}$.
For $\varphi\in\X$ the function $\varphi_{k,j}=\varphi\chi_k\theta_j(\eta_1)$ lies in $W^{1,2}(\Rp)$ and has compact support in $\Rp$, so (a) gives
$\mathcal B(h;\varphi_{k,j})=0$. Now $\nabla\varphi_{k,j}=\chi_k\theta_j\nabla\varphi+\varphi\theta_j\nabla\chi_k+\varphi\chi_k\theta_j'e_1$, and we let $j,k\to\infty$ in each term of $\mathcal B(h;\varphi_{k,j})=0$.
\begin{itemize}[leftmargin=*]
\item $\int\nabla h\cdot\nabla\varphi\,\chi_k\theta_jK\to\int\nabla h\cdot\nabla\varphi\,K$ by dominated convergence ($\nabla h,\nabla\varphi\in L^2_K$).
\item $\big|\int\nabla h\cdot\nabla\chi_k\,\varphi\theta_jK\big|\le2\|\nabla h\|_{L^2_K}\|\varphi\|_{L^2_K(\{r>k\})}\to0$.
\item For the strip term we use Hardy's inequality. By the identity $\int|\nabla\phi|^2K=\int|\nabla v|^2+(\frac{|\eta|^2}{16}+\frac N4)v^2$, $v=K^{1/2}\phi$, valid for $\phi\in C_c^\infty(\Rp)$ (proof of Lemma~\ref{lem:reg}), the map $\phi\mapsto K^{1/2}\phi$ sends $W^{1,2}_{0,K}(\Rp)$ into $W^{1,2}_0(\Rp)$. Hence $v=K^{1/2}\varphi\in W^{1,2}_0(\Rp)$ and $\int v^2/\eta_1^2\le4\int|\partial_1v|^2<\infty$ (one-dimensional Hardy on lines). Since $j<2/\eta_1$ on $\{\frac1j<\eta_1<\frac2j\}$,
\[
\Big|\int\partial_1h\,\varphi\chi_k\theta_j'K\Big|\le\|\nabla h\|_{L^2_K}\Big(j^2\!\!\int_{\{\frac1j<\eta_1<\frac2j\}}\!\!\varphi^2K\Big)^{1/2}\le2\|\nabla h\|_{L^2_K}\Big(\int_{\{\eta_1<2/j\}}\frac{v^2}{\eta_1^2}\Big)^{1/2}\to0,
\]
uniformly in $k$.
\item $\int\eta_1^{2\alpha}h^q\varphi_{k,j}K\to\int\eta_1^{2\alpha}h^q\varphi K$ and $\lambda\int h\varphi_{k,j}K\to\lambda\int h\varphi K$ by dominated convergence, because $\eta_1^{2\alpha}h^q|\varphi|\in L^1(K)$ by H\"older ($h,\varphi\in L^{q+1}_{\eta_1^{2\alpha}K}$) and $h\varphi\in L^1(K)$.
\end{itemize}
This gives \eqref{eq:weak}.
\end{proof}

%---------------------------------------------------------------------
\section{Explicit barriers}

For $w=e^{-r^2/4}g(r)\psi(\sigma)$ a direct computation gives
\begin{equation}\label{eq:radang}
L_Kw-\lambda w=-e^{-r^2/4}\Big[\mathcal Dg\cdot\psi+r^{-2}g\,\Delta_{\Sph}\psi\Big],\qquad
\mathcal Dg=g''+\Big(\frac{N-1}r-\frac r2\Big)g'+\Big(\lambda-\frac N2\Big)g,
\end{equation}
and for $g=r^p$,
\begin{equation}\label{eq:Dp}
\mathcal Dr^p=\tfrac12(s-p)\,r^p+p(p+N-2)\,r^{p-2}.
\end{equation}
For a function of $\sigma_1$ alone, $\Delta_{\Sph}f(\sigma_1)=(1-\sigma_1^2)f''-(N-1)\sigma_1f'$;
in particular $\Delta_{\Sph}\sigma_1=-(N-1)\sigma_1$.

\begin{lemma}[exact solution]\label{lem:W}
Let $\beta=2\lambda$ and $M^{q-1}=\beta(\beta+1)$. Then $W(\eta)=M\eta_1^{-\beta}$ satisfies
$\mathcal M(W)=0$ in $\Rp$.
\end{lemma}
\begin{proof}
$-W''=-M\beta(\beta+1)\eta_1^{-\beta-2}$; $-\frac{\eta_1}2W'-\lambda W=(\frac\beta2-\lambda)W=0$;
and $2\alpha-\beta q=-\beta-2$ because $\beta(q-1)=2\alpha+2$, so
$\eta_1^{2\alpha}W^q=M^q\eta_1^{-\beta-2}$. Therefore $\mathcal M(W)=M\eta_1^{-\beta-2}(M^{q-1}-\beta(\beta+1))=0$.
\end{proof}
$W$ is the self-similar form of the stationary solution $Mx_1^{-2(1+\alpha)/(q-1)}$.

\begin{remark}\label{rem:nosuper}
The function $U=z_1(r)\sigma_1$, with $z_1$ the recessive solution of the radial equation
$z''+(\frac{N-1}r+\frac r2)z'+(\lambda-\frac{N-1}{r^2})z=0$, cannot be positive on $(0,\infty)$ if
$\lambda>\lambda_1$. Indeed let $w=re^{-r^2/4}$ (the radial part of $\psi_1$) and $p=r^{N-1}e^{r^2/4}$.
Then $\big(p(z_1'w-z_1w')\big)'=-(\lambda-\lambda_1)pz_1w$. The boundary term tends to $0$ at $\infty$,
and to $-Nb$ at $0$ if $z_1\sim ar+br^{-(N-1)}$. If $z_1>0$ then $b\ge0$, while integration gives
$Nb=-(\lambda-\lambda_1)\int pz_1w<0$, a contradiction. So $z_1\sigma_1$ is not a global supersolution, which is why we use the explicit solution $W$ of Lemma~\ref{lem:W} instead.
Marcus and V\'eron \cite[proof of Lemma~2.5, Steps~2--3]{MV01} use $z_1\sigma_1$ only in $\{r>R_0\}$, where $z_1>0$, which is correct.
\end{remark}

\begin{lemma}[upper barrier]\label{lem:U}
Let $A_0=|(s-1)(s+N-1)|+1$. There is $R_0\ge\sqrt{2A_0}$ such that
\[
U(\eta)=\sigma_1\,r^{s}\,(1-A_0r^{-2})\,e^{-r^2/4}
\]
satisfies $\mathcal M(U)\ge0$ in $\{r>R_0\}$, and
$\tfrac12\sigma_1r^se^{-r^2/4}\le U\le\sigma_1r^se^{-r^2/4}$ there.
\end{lemma}
$U$ is an explicit substitute for the supersolution $z_1(r)\sigma_1$ used, for $\alpha=0$, in \cite[proof of Lemma~2.5]{MV01}, where $z_1$ is defined through an ODE. The lower barrier below is new; it is what gives the lower bound in \eqref{eq:sharp} with the exponent of the upper bound.
\begin{proof}
By \eqref{eq:radang}--\eqref{eq:Dp} with $\psi=\sigma_1$, $g=r^s-A_0r^{s-2}$,
\[
L_KU-\lambda U=e^{-r^2/4}\sigma_1\Big[\big(A_0-(s-1)(s+N-1)\big)r^{s-2}+A_0c_3\,r^{s-4}\Big],
\]
with $c_3=(s-2)(s+N-4)-(N-1)$. Here we used $s(s+N-2)-(N-1)=(s-1)(s+N-1)$. The first coefficient is
$\ge1$, so the bracket is $\ge\frac12r^{s-2}$ for $r\ge R_0$ large. The term
$\eta_1^{2\alpha}U^q$ is nonnegative.
\end{proof}

\begin{lemma}[lower barrier]\label{lem:w}
Assume $1<q<q_c$ (so $s>1$) and put $m=\min\{2,\,2\alpha+q+2\}$; then $m>1$. There exist
$\tau\in(0,1]$ and $R_1\ge1$ such that for every $\varepsilon\in(0,1]$ the function
\[
w(\eta)=\varepsilon\,\psi(\sigma)\,r^{s}e^{-r^2/4},\qquad\psi=\sigma_1+\tau\sigma_1^{m},
\]
satisfies $\mathcal M(w)\le0$ in $\{r>R_1\}$. Moreover $\sigma_1\le\psi\le(1+\tau)\sigma_1$.
\end{lemma}
\begin{proof}
By \eqref{eq:radang}--\eqref{eq:Dp} with $p=s$, $L_Kw-\lambda w=-\varepsilon e^{-r^2/4}r^{s-2}G(\sigma_1)$, where
$G=s(s+N-2)\psi+\Delta_{\Sph}\psi$. With $t=\sigma_1$,
\[
G(t)=(s-1)(s+N-1)\,t+\tau\Big[m(m-1)t^{m-2}(1-t^2)+\big(s(s+N-2)-m(N-1)\big)t^m\Big].
\]
Put $A_1=(s-1)(s+N-1)>0$, $B_1=m(m-1)>0$, $D_1=|s(s+N-2)-m(N-1)|$. Choose $t_0\in(0,1)$ with
$D_1t_0^2\le\frac12B_1(1-t_0^2)$, then $\tau\in(0,1]$ with $\tau D_1\le\frac12A_1t_0$. For
$t\le t_0$: $G\ge\tau t^{m-2}\big[B_1(1-t_0^2)-D_1t_0^2\big]\ge\frac12\tau B_1(1-t_0^2)t^{m-2}$.
For $t\in[t_0,1]$: $G\ge A_1t_0-\tau D_1\ge\frac12A_1t_0\ge\frac12A_1t_0^{3-m}t^{m-2}$. Hence
$G(t)\ge c_2t^{m-2}$ on $(0,1]$ for some $c_2>0$. Next,
$\eta_1^{2\alpha}w^q\le\varepsilon^q(1+\tau)^q r^{2\alpha+sq}\sigma_1^{2\alpha+q}e^{-qr^2/4}$ and
$\sigma_1^{2\alpha+q}\le\sigma_1^{m-2}$ because $m-2\le2\alpha+q$. Therefore
\[
\mathcal M(w)\le\varepsilon e^{-r^2/4}r^{s-2}\sigma_1^{m-2}\Big[-c_2+\varepsilon^{q-1}(1+\tau)^q\,r^{2\alpha+s(q-1)+2}e^{-(q-1)r^2/4}\Big]\le0
\]
for $r\ge R_1$, uniformly in $\sigma$ and in $\varepsilon\le1$.
\end{proof}
Note where subcriticality enters: the coefficient $A_1=(s-1)(s+N-1)$ is positive exactly when $s>1$, i.e.\ $q<q_c$.

\begin{remark}[a global supersolution in the energy space]\label{rem:g}
The exact solution $W$ blows up at $\partial\Rp$, which is why existence is proved by exhaustion in Section~4. A global supersolution that belongs to $\X$ is also available. Fix $p\in(\frac12,1)$ and $k>2\lambda-N-p$, and let
\[
g(\eta)=\eta_1^{p}\,(1+|\eta|^2)^{k/2}\,e^{-|\eta|^2/4}.
\]
A direct computation gives, with $r=|\eta|$,
\[
\frac{(L_K-\lambda)g}{g}=\frac{p(1-p)}{\eta_1^2}-\frac{k(N+2p)}{1+r^2}-\frac{k(k-2)r^2}{(1+r^2)^2}+\frac p2+\frac{kr^2}{2(1+r^2)}+\frac N2-\lambda .
\]
The right side tends to $+\infty$ as $\eta_1\to0$ (locally uniformly in $\eta'$) and to $\frac{N+p+k}2-\lambda>0$ as $r\to\infty$, so it is negative only on a compact set $Q\Subset\Rp$. Since $\eta_1^{2\alpha}g^{q-1}$ has a positive minimum on $Q$, the function $Cg$ satisfies $\mathcal M(Cg)\ge0$ in $\Rp$ for $C$ large. Moreover $g\in\X$: near $\partial\Rp$, $|\nabla g|^2\sim\eta_1^{2p-2}$ is integrable because $p>\frac12$, and $\eta_1^{2\alpha}g^{q+1}\sim\eta_1^{2\alpha+p(q+1)}$ with $2\alpha+p(q+1)>-1$. So $Cg$ can replace $W$ as the upper barrier in the construction of Section~4, with subsolution $\varepsilon\phi_D$, where $\phi_D$ is the first Dirichlet eigenfunction of $L_K$ in a large bounded domain $D\Subset\Rp$ (its eigenvalue is $<\lambda$ by Lemma~\ref{lem:spec}(i)). This gives the bound $h\le Cg$, which has the Gaussian factor but an algebraic exponent larger than the sharp one in \eqref{eq:sharp}.
\end{remark}

%---------------------------------------------------------------------
\section{Existence for $1<q<q_c$}

From now on $1<q<q_c$, so $\lambda>\lambda_1$. Let $D_n=\Rp\cap B_n(0)$.

The construction by exhaustion follows \cite[proof of Lemma~2.5, Step~2]{MV01}, where it is used for $\alpha=0$ in exterior domains; the new element is the global upper bound by the exact solution $W$ of Lemma~\ref{lem:W}.

\begin{proposition}[approximate solutions]\label{prop:hn}
There is $n_0$ such that for $n\ge n_0$ there exists $h_n\in W^{1,2}_0(D_n)\cap L^\infty(D_n)$,
positive in $D_n$, satisfying \eqref{eq:weak} for all
$\varphi\in X_n:=\{\varphi\in W_0^{1,2}(D_n):\int\eta_1^{2\alpha}|\varphi|^{q+1}<\infty\}$.
Its zero extension $\hat h_n$ is a weak subsolution in $D_{n+1}$.
\end{proposition}
\begin{proof}
Minimize $J(h)=\frac12\int|\nabla h|^2K-\frac\lambda2\int h^2K+\frac1{q+1}\int\eta_1^{2\alpha}|h|^{q+1}K$ on $X_n$.

\emph{Coercivity.} Fix $\delta=\delta(n)$ small. On $\{\eta_1<\delta\}$ the one-dimensional Poincar\'e inequality
($h=0$ at $\eta_1=0$) gives $\int h^2\le\frac{4\delta^2}{\pi^2}\int|\partial_1h|^2$, and since
$K\le e^{n^2/4}$ on $D_n$, we get $\frac\lambda2\int_{\{\eta_1<\delta\}}h^2K\le\frac14\int|\nabla h|^2K$.
On $\{\eta_1\ge\delta\}\cap D_n$ the weight $\eta_1^{2\alpha}$ is bounded below by a positive constant, so Young's
inequality gives $\frac\lambda2h^2\le\frac1{2(q+1)}\eta_1^{2\alpha}|h|^{q+1}+C$. Hence
$J(h)\ge\frac14\int|\nabla h|^2K+\frac1{2(q+1)}\int\eta_1^{2\alpha}|h|^{q+1}K-C_n$. No integrability of
$\eta_1^{2\alpha}$ near $\partial\Rp$ is needed. Rellich's theorem and convexity of the other terms give a minimizer.

\emph{Nontriviality.} By Lemma~\ref{lem:spec}, $\lambda_1(D_n)<\lambda$ for $n\ge n_0$. Take $\phi\in C_c^\infty(D_n)$ with Rayleigh
quotient $<\lambda$; then $J(t\phi)<0$ for small $t$.

\emph{Equation.} Since $J(|h|)=J(h)$ we may take $h\ge0$. For $\varphi\in X_n$ the map $t\mapsto J(h+t\varphi)$ is differentiable by
dominated convergence, which gives \eqref{eq:weak} on $X_n$. Testing with $(h-k)_+$ and Stampacchia's method (as in \cite[Theorem~8.15]{GT}) give
$h\in L^\infty$. The strong maximum principle \cite[Theorem~3.5]{GT} gives $h>0$ in $D_n$. \emph{Subsolution.} For $0\le\varphi\in W_0^{1,2}(D_{n+1})\cap L^\infty$ test the equation of $h_n$ with $\varphi\,T_\tau(h_n)$,
$T_\tau(t)=\min(t/\tau,1)$, drop the nonnegative term $\int\varphi T_\tau'(h_n)|\nabla h_n|^2K$, and let $\tau\to0$ (dominated convergence for the linear terms, monotone convergence for the absorption term, whose limit is therefore finite).
\end{proof}

\begin{proposition}[bounds uniform in $n$]\label{prop:unif}
There are constants independent of $n\ge n_0$ such that
\begin{enumerate}[label=\textup{(\roman*)},leftmargin=*]
\item $h_n\le h_{n+1}$;
\item $h_n\le W$ in $D_n$;
\item $h_n\le C_1\min(\eta_1,1)$ in $D_n$;
\item $h_n\le C_2U$ in $\Rp\cap\{R_0<r<n\}$;
\item $\int|\nabla h_n|^2K+\int\eta_1^{2\alpha}h_n^{q+1}K\le C_3$.
\end{enumerate}
\end{proposition}
\begin{proof}
Since $h_n$ is only known to lie in $W_0^{1,2}(D_n)$, hypothesis (H3) is checked below through the zero trace of $h_n$ on $\partial D_n$, and,
at the interior parts of $\partial Q$ (the faces $\{\eta_1=\delta\}$ and $\{r=R_0\}$), through pointwise inequalities valid in a neighbourhood.

(i) Lemma~\ref{lem:comp} in $Q=D_{n+1}$ with $u_1=\hat h_n$, $u_2=h_{n+1}$. Indeed $(\hat h_n-h_{n+1}-\varepsilon)_+\le(\hat h_n-\varepsilon)_+\in W_0^{1,2}(D_{n+1})$.

(ii) Lemma~\ref{lem:comp} in $Q=D_n$ with $u_1=h_n$, $u_2=W$. Here $S_\varepsilon\subset\{W<\|h_n\|_\infty\}$ stays away from $\partial\Rp$, and $W>0$.

(iii) Let $\delta=\pi/(2\sqrt{\lambda+1/2})$, $\kappa=\pi/(2\delta)$ (so $\kappa^2=\lambda+\frac12$), and $\Phi(\eta_1)=M\delta^{-\beta}\sin(\kappa\eta_1)$.
Since $x\cos x\le\sin x$ on $[0,\pi/2]$, $-\frac{\eta_1}2\Phi'\ge-\frac12\Phi$, so
$L_K\Phi-\lambda\Phi\ge(\kappa^2-\frac12-\lambda)\Phi\ge0$ and $\mathcal M(\Phi)\ge0$ in $Q=D_n\cap\{\eta_1<\delta\}$.
On $\{\eta_1=\delta\}$, (ii) gives $h_n\le M\delta^{-\beta}=\Phi(\delta)$; on $\{\eta_1=0\}$ and on the sphere $h_n=0$.
Lemma~\ref{lem:comp} gives $h_n\le\Phi\le M\delta^{-\beta}\min\{\kappa\eta_1,1\}$ for $\eta_1<\delta$. For $\eta_1\ge\delta$,
(ii) gives $h_n\le M\delta^{-\beta}$. With $C_1=M\delta^{-\beta}\max\{\kappa,1/\delta,1\}$ this gives (iii) (the case $\delta>1$, $1<\eta_1<\delta$, needs the $\max$).

(iv) Choose $C_2$ so large that $C_2U\ge C_1r\sigma_1\ge h_n$ on $\{R_0\le r\le R_0+1\}$ (possible by Lemma~\ref{lem:U} and (iii)).
Apply Lemma~\ref{lem:comp} in $\Rp\cap\{R_0<r<n\}$ with $u_1=h_n$, $u_2=C_2U$; $C_2U$ is a supersolution since $\mathcal M(C_2U)\ge C_2(\ldots)\ge0$.

(v) Test the equation of $h_n$ with $h_n$ itself: $\int|\nabla h_n|^2K+\int\eta_1^{2\alpha}h_n^{q+1}K=\lambda\int h_n^2K$. By (iii) and (iv),
$h_n\le\bar H:=C_1\min(\eta_1,1)\mathbf 1_{\{r<R_0\}}+C_2U\mathbf 1_{\{r\ge R_0\}}\in L^2_K(\Rp)$,
because $U^2K\lesssim r^{2s}e^{-r^2/4}$.
\end{proof}

\begin{proposition}[the solution]\label{prop:h}
$h=\lim h_n$ exists and $h\in\Hh\cap C^2(\Rp)\cap C(\overline{\Rp})$, $h>0$ in $\Rp$, $h=0$ on $\partial\Rp$, and
$h\le C_1\min(\eta_1,1)$, $h\le C_2U$ for $r\ge R_0$.
\end{proposition}
\begin{proof}
By (i) and (ii) the sequence is increasing and locally bounded in $\Rp$. Interior estimates give $C^2_{\rm loc}(\Rp)$ convergence and
that $h$ solves \eqref{eq:E} in $\Rp$. Continuity at $\partial\Rp$ follows from $h\le C_1\eta_1$. Positivity follows from $h\ge h_{n_0}>0$.
By (v) and (iii)--(iv), $(h_n)$ is bounded in $W^{1,2}_{0,K}(\Rp)$ (extended by $0$) and $h_n\to h$ pointwise with $h_n\le h\le\bar H\in L^2_K$; hence the whole sequence converges weakly in $W^{1,2}_{0,K}(\Rp)$ to $h$,
so $h\in W^{1,2}_{0,K}(\Rp)$, and $h\in L^{q+1}_{\eta_1^{2\alpha}K}$ by Fatou and (v). Thus $h\in\X$, and $h\to0$ at infinity by (iv).

\emph{The identity \eqref{eq:weak}.} Let $\varphi\in\X$ have support in $B_\rho$. Take $\varphi_i\in C_c^\infty(\Rp)$ with $\varphi_i\to\varphi$ in $W^{1,2}_K$ and a cut-off $\chi$ ($\chi=1$ on $B_\rho$, supported in $B_{\rho+1}$); then $\chi\varphi_i\in C_c^\infty(D_n)$ for $n>\rho+1$ and $\chi\varphi_i\to\varphi$ in $W^{1,2}$, so $\varphi\in W^{1,2}_0(D_n)$.
Since $\int\eta_1^{2\alpha}|\varphi|^{q+1}<\infty$, $\varphi\in X_n$, and Proposition~\ref{prop:hn} gives
$\int\nabla h_n\cdot\nabla\varphi K+\int\eta_1^{2\alpha}h_n^q\varphi K=\lambda\int h_n\varphi K$. Let $n\to\infty$: the first term converges by weak convergence in $L^2_K$ of $\nabla h_n$; the other two converge by dominated convergence, since $0\le h_n\le h$ and
$\eta_1^{2\alpha}h^q|\varphi|\in L^1(K)$ (H\"older), $h\varphi\in L^1(K)$. For general $\varphi\in\X$ apply this to $\varphi\chi_k(r)\in\X$ (bounded support) and let $k\to\infty$: $\nabla(\varphi\chi_k)=\chi_k\nabla\varphi+\varphi\nabla\chi_k$, and
$\big|\int\nabla h\cdot\nabla\chi_k\,\varphi K\big|\le2\|\nabla h\|_{L^2_K}\|\varphi\|_{L^2_K(\{r>k\})}\to0$, while the remaining terms pass to the limit by dominated convergence. Hence \eqref{eq:weak} holds for all $\varphi\in\X$, i.e.\ $h\in\Hh$.
\end{proof}

\begin{lemma}[linear lower bound]\label{lem:hopf}
For every $R>0$ there is $c_R>0$ with $h(\eta)\ge c_R\,\eta_1$ for $\eta\in\Rp$, $|\eta|\le R$.
\end{lemma}
\begin{proof}
Let $\gamma=2(N-1)+R+1$, $\mu=2\gamma(N+R)/(e^{-\gamma/4}-e^{-\gamma})$, $m=\min\{2,2\alpha+q+2\}\in(1,2]$, $\rho=\min\{1,(m(m-1)/(4\mu))^{1/(3-m)}\}$ and
$\varepsilon_0=\min\{1,(m(m-1)/2^{q+1})^{1/(q-1)}\}$. For $|z|\le1$, $z\in\R^{N-1}$, let $\phi(z)=(e^{-\gamma|z|^2}-e^{-\gamma})/(1-e^{-\gamma})$, so $0\le\phi\le1$, $\phi(0)=1$, $\phi=0$ on $|z|=1$.
Fix $y'$ with $|y'|\le R$, put $z=\eta'-y'$ and $L'=-\Delta'-\frac12\eta'\cdot\nabla'$. Then
$L'\phi=(1-e^{-\gamma})^{-1}e^{-\gamma|z|^2}[2\gamma(N-1)-4\gamma^2|z|^2+\gamma\eta'\cdot z]$ with $|\eta'\cdot z|\le(R+1)|z|$; the bracket is $\le0$ for $|z|\ge\frac12$ and $\le2\gamma(N+R)$ for $|z|\le\frac12$,
where $\phi\ge(e^{-\gamma/4}-e^{-\gamma})/(1-e^{-\gamma})$. Hence $L'\phi\le\mu\phi$ on $\{|z|<1\}$.

Let $Q=\{0<\eta_1<\rho,\ |z|<1\}$ and $\Theta=\theta(\eta_1)\phi(z)$, $\theta(t)=\varepsilon(t+t^m)$, $\varepsilon\le\varepsilon_0$. Since $K$ separates, $L_K(\theta\phi)=\phi L_1\theta+\theta L'\phi$ with $L_1\theta=-\theta''-\frac{\eta_1}2\theta'$.
Using $\theta\ge0$, $L'\phi\le\mu\phi$, $\phi^q\le\phi$, $-\lambda\theta\le0$ and $\theta'>0$, for $\eta_1\le\rho\le1$,
\[
\mathcal M(\Theta)\le\phi\big[-\theta''+\mu\theta+\eta_1^{2\alpha}\theta^q\big]\le\varepsilon\phi\,\eta_1^{m-2}\big[-m(m-1)+2\mu\eta_1^{3-m}+2^q\varepsilon^{q-1}\eta_1^{2\alpha+q+2-m}\big]\le0,
\]
because $3-m\ge1$ and $2\alpha+q+2-m\ge0$, and by the choice of $\rho$ and $\varepsilon_0$. Let $m_R=\min\{h:\eta_1=\rho,\ |\eta'|\le R+1\}>0$ and $\varepsilon=\min\{\varepsilon_0,m_R/(2\rho)\}$. Then
$\Theta\le2\varepsilon\rho\le h$ on $\{\eta_1=\rho\}$, and $\Theta=0$ on the other faces of $Q$. Since $\Theta\le2\varepsilon\eta_1$ and $\Theta-h$ is continuous and $\le0$ on the top and lateral faces, each set $\{\Theta>h+\varepsilon'\}$ has compact closure in $Q$;
Lemma~\ref{lem:comp} (with Lemma~\ref{lem:test}(a)) gives $\Theta\le h$ in $Q$, so $h(\eta_1,y')\ge\varepsilon\eta_1$ for $\eta_1<\rho$. For $\eta\in\Rp$, $|\eta|\le R$: if $\eta_1<\rho$ take $y'=\eta'$; if $\eta_1\ge\rho$ use $h\ge\min\{h:\eta_1\ge\rho,\ |\eta|\le R\}>0$ and $\eta_1\le R$.
\end{proof}

\begin{proposition}[sharp lower bound]\label{prop:lower}
There are $c,R_1>0$ with $h\ge c\,\sigma_1r^{s}e^{-r^2/4}$ for $r\ge R_1$.
\end{proposition}
\begin{proof}
Take $\tau,R_1$ from Lemma~\ref{lem:w}. By Lemma~\ref{lem:hopf}, $h\ge c_{R_1}R_1\sigma_1$ on $\{r=R_1\}$. Choose $\varepsilon\le1$ with
$(1+\tau)\varepsilon R_1^se^{-R_1^2/4}\le c_{R_1}R_1$, so $w\le h$ on $\{r=R_1\}$. Apply Lemma~\ref{lem:comp} in $Q=\Rp\cap\{r>R_1\}$ with
$u_1=w$ and $u_2=h$: $w$ is a subsolution, $w\in L^2_K\cap L^{q+1}_{\eta_1^{2\alpha}K}$ (use $2\alpha+q+1>0$), and since $w,h$ are
continuous on $\overline{\Rp}$, vanish on $\partial\Rp$ and $w\le h$ on $\{r=R_1\}$, the sets $S_\varepsilon$ are bounded and stay away from
$\partial\Rp\cup\{r=R_1\}$.
\end{proof}

%---------------------------------------------------------------------
\section{Uniqueness and nonexistence}

For $\alpha=0$, Propositions~\ref{prop:uniq} and~\ref{prop:nonex} are \cite[proof of Lemma~2.5, Step~1]{MV01} and \cite[proof of Theorem~2.1(i)]{MV01}; there uniqueness is proved by testing the equations against each other, here by Lemma~\ref{lem:comp}, and nonexistence by the same comparison with $\lambda_1$.

\begin{proposition}[uniqueness]\label{prop:uniq}
Let $1<q<q_c$ and let $h'\in\Hh$, $h'\not\equiv0$. Then $h'=h$.
\end{proposition}
\begin{proof}
$h'$ is bounded (Lemma~\ref{lem:reg}), hence $C^{2,\theta}_{\rm loc}(\Rp)$ by interior elliptic regularity \cite[Theorem~8.8, Lemma~9.16, Theorem~9.19]{GT}, and $h'>0$ in $\Rp$ by the strong maximum principle \cite[Theorem~3.5]{GT}. For $n\ge n_0$, Lemma~\ref{lem:comp} in $D_n$ with $u_1=h_n$, $u_2=h'$ gives $h_n\le h'$, so
$h\le h'$. Now apply Lemma~\ref{lem:comp} in $Q=\Rp$ with $u_1=h'$, $u_2=h$: (H2) holds because $h'\in\X$, and $h'$ is bounded by Lemma~\ref{lem:reg}; in (H3), $S_\varepsilon$ is bounded because
$h'\to0$ at infinity (Lemma~\ref{lem:reg}), and $(h'-h-\varepsilon)_+\le(h'-\varepsilon)_+\in W_0^{1,2}$. Hence $h'\le h$.
\end{proof}

\begin{proposition}[nonexistence]\label{prop:nonex}
Let $q\ge q_c$, i.e.\ $\lambda\le\lambda_1$. If $h\in\Hh$ then $h\equiv0$.
\end{proposition}
\begin{proof}
$\psi_1\in\X$, so \eqref{eq:weak} holds with $\varphi=\psi_1$. Since $\psi_1$ is an eigenfunction and $h\in W^{1,2}_{0,K}$,
$\int\nabla h\cdot\nabla\psi_1K=\lambda_1\int h\psi_1K$. Hence
\[
(\lambda-\lambda_1)\int_Ah\psi_1K=\int_{\Rp}\eta_1^{2\alpha}h^q\psi_1K\ge0 .
\]
If $\lambda\le\lambda_1$ the left side is $\le0$, so both sides vanish, and $h^q\psi_1=0$ a.e., i.e.\ $h\equiv0$.
\end{proof}

%---------------------------------------------------------------------
\section{Continuous profiles: uniqueness among all self-similar solutions}\label{sec:cont}

Let $\mathcal S$ be the set of $h\in C^2(\Rp)\cap C(\overline{\Rp})$ with $h>0$ in $\Rp$, solving the equation in \eqref{eq:E} classically, with $h=0$ on $\partial\Rp$ and
\begin{equation}\label{eq:decay}
\lim_{|\eta|\to\infty}|\eta|^{2\lambda}h(\eta)=0 .
\end{equation}
If $u\in C(\overline E\setminus\{(0,0)\})$ is a positive self-similar solution of \eqref{eq:pde}--\eqref{eq:bc}, i.e.\ $u=t^{-\lambda}h(x/\sqrt t)$, then
$h=u(\cdot,1)\in\mathcal S$: only \eqref{eq:decay} needs a comment, and it follows from $|\eta|^{2\lambda}h(\eta)=u(\eta/|\eta|,|\eta|^{-2})$ and the uniform
continuity of $u$ on the compact set $\{|x|=1,\ x_1\ge0,\ 0\le t\le1\}$, where $u(x,0)=0$. Conversely every $h\in\mathcal S$ gives such a solution.
The point of this section is that no energy assumption is needed: $\mathcal S\subset\Hh$ for every $q>1$. For $\alpha=0$ this is \cite[proof of Lemma~2.5, Steps~3--4]{MV01}; the weight requires the linear boundary bound of Lemma~\ref{lem:bdry} and the Hardy argument in Step~3 below.

\begin{lemma}[linear boundary bound]\label{lem:bdry}
Let $h\in C^2(\Rp)\cap C(\overline{\Rp})$, $h\ge0$, with $\mathcal M(h)\le0$ in $\Rp$ and $h=0$ on $\partial\Rp$. For every $R>0$ there is $C_R$ such that
$h(\eta)\le C_R\,\eta_1$ for all $\eta\in\Rp$ with $|\eta|\le R$.
\end{lemma}
\begin{proof}
Let $m=\sup\{h(\eta):\eta\in\Rp,\ |\eta|\le R+2\}<\infty$ (continuity on a compact set); we may assume $m>0$. Fix $y=(0,y')\in\partial\Rp$ with $|y|\le R$ and write $z=\eta'-y'$.
For $\rho\in(0,(2N)^{-1/2}]$ put $\rho_2=\sqrt{2N}\,\rho\le1$ and
\[
Q=\{0<\eta_1<\rho,\ |z|<\rho_2\}\subset\Rp\cap B_{R+2},\qquad
\Theta(\eta)=a\eta_1-b\eta_1^2+c|z|^2,\quad a=\frac{2m}{\rho},\ b=\frac m{\rho^2},\ c=\frac m{\rho_2^2}.
\]
\emph{$\Theta$ is a positive supersolution.} Since $a-b\eta_1\ge a/2$ on $Q$, $\Theta>0$ in $Q$. Moreover $\Delta\Theta=-2b+2c(N-1)$ and $\eta\cdot\nabla\Theta=\eta_1(a-2b\eta_1)+2c\,\eta'\!\cdot z$, so
\[
L_K\Theta=2b-2c(N-1)-\tfrac12a\eta_1+b\eta_1^2-c\,\eta'\!\cdot z .
\]
On $Q$ we have $0<\eta_1<\rho$, $|z|<\rho_2$, $|\eta'|\le R+1$; dropping the nonnegative terms $b\eta_1^2$ and $\lambda b\eta_1^2$,
\[
\begin{aligned}
(L_K-\lambda)\Theta&\ge\ 2b-2c(N-1)-\tfrac12a\rho-c(R+1)\rho_2-\lambda(a\rho+c\rho_2^2)\\
&=\frac m{\rho^2}\Big(2-\frac{N-1}{N}\Big)-m\Big(1+3\lambda+\frac{R+1}{\sqrt{2N}\,\rho}\Big)\ \ge\ \frac m{2\rho^2},
\end{aligned}
\]
where the last inequality holds for $\rho=\rho_*(N,R,\lambda)$ small enough, because $2-\frac{N-1}N\ge1$. Hence $\mathcal M(\Theta)=(L_K-\lambda)\Theta+\eta_1^{2\alpha}\Theta^q\ge0$ in $Q$.

\emph{Boundary values.} On $\{\eta_1=0\}$: $\Theta=c|z|^2\ge0=h$. On $\{|z|=\rho_2\}$: $\Theta\ge c\rho_2^2=m\ge h$ (because $a\eta_1-b\eta_1^2\ge0$). On $\{\eta_1=\rho\}$: $\Theta\ge a\rho-b\rho^2=m\ge h$. So $h\le\Theta$ on $\partial Q$.

\emph{Comparison.} Apply Lemma~\ref{lem:comp} in $Q$ with $u_1=h$, $u_2=\Theta$. Since $h,\Theta$ are continuous on $\overline Q$ and $h\le\Theta$ on $\partial Q$, each $S_\varepsilon=\{h>\Theta+\varepsilon\}$ is relatively compact in $Q$; hence (H1)--(H3) hold with classical sub- and supersolutions, and (H2) is automatic. Therefore $h\le\Theta$ in $Q$.

\emph{Conclusion.} On the segment $z=0$: $h(\eta_1,y')\le a\eta_1=(2m/\rho_*)\eta_1$ for $\eta_1<\rho_*$. Given $\eta\in\Rp$ with $|\eta|\le R$ and $\eta_1<\rho_*$, take $y=(0,\eta')$. If $\eta_1\ge\rho_*$, then $h\le m\le(m/\rho_*)\eta_1$. So $C_R=2m/\rho_*$ works.
\end{proof}

\begin{lemma}[decaying supersolution]\label{lem:Z}
Let $b=\max\{0,2\lambda(2\lambda+2-N)\}+1$ and $Z=r^{-2\lambda}+b\,r^{-2\lambda-2}$. There is $R_0'$ such that $L_KZ-\lambda Z>0$, hence $\mathcal M(Z)>0$, in $\{r\ge R_0'\}$.
\end{lemma}
\begin{proof}
For a power, $L_Kr^{-p}-\lambda r^{-p}=(\frac p2-\lambda)r^{-p}+p(N-2-p)r^{-p-2}$. With $p=2\lambda$ the first term vanishes; with $p=2\lambda+2$ it equals $r^{-2\lambda-2}$. So
\[
L_KZ-\lambda Z=\big[b-2\lambda(2\lambda+2-N)\big]r^{-2\lambda-2}+b(2\lambda+2)(N-4-2\lambda)r^{-2\lambda-4}\ \ge\ r^{-2\lambda-2}\big(1-C r^{-2}\big),
\]
which is positive for $r$ large.
\end{proof}

\begin{proposition}\label{prop:S}
Let $q>1$. Every $h\in\mathcal S$ satisfies: $h\le C_R\eta_1$ for $|\eta|\le R$; $h\le CU$ for $r\ge R_0$ (with $U$ from Lemma~\ref{lem:U}); and $h\in\Hh$.
\end{proposition}
\begin{proof}
\emph{Step 1.} $h$ is a nonnegative solution, so Lemma~\ref{lem:bdry} gives $h\le C_R\eta_1$ on $\{|\eta|\le R\}$ for every $R$.

\emph{Step 2.} Take $R_0$ so large that Lemmas~\ref{lem:U} and~\ref{lem:Z} hold in $\{r\ge R_0\}$, and choose $C$ with $CU\ge h$ on $\{R_0\le r\le R_0+1\}$ (Step 1 and $U\ge\frac12\sigma_1r^se^{-r^2/4}$).
Given $\varepsilon>0$ use \eqref{eq:decay} to pick $R_\varepsilon>R_0+1$ with $h\le\varepsilon r^{-2\lambda}\le\varepsilon Z$ for $r\ge R_\varepsilon$. Then $V_\varepsilon=CU+\varepsilon Z>0$ satisfies $\mathcal M(V_\varepsilon)\ge0$ in $\{r>R_0\}$,
and $V_\varepsilon\ge h$ on $\{R_0\le r\le R_0+1\}\cup\{r\ge R_\varepsilon\}$. Apply Lemma~\ref{lem:comp} in $Q=\Rp\cap\{R_0<r<R_\varepsilon+1\}$ with $u_1=h$, $u_2=V_\varepsilon$: the sets $S_{\varepsilon'}$ stay in $\{R_0+1<r<R_\varepsilon\}$ and away from $\partial\Rp$
(where $h\to0$ while $V_\varepsilon\ge\varepsilon Z>0$ is bounded below), so they are relatively compact in $Q$ (and in $\Rp$), which gives (H1)--(H3). Hence $h\le V_\varepsilon$ in $Q$. Letting $\varepsilon\to0$ gives $h\le CU$ in $\{r>R_0\}$.

\emph{Step 3 (energy).} By Steps 1--2, $h$ is bounded and $\int_Ah^2K<\infty$ (use $U^2K\lesssim r^{2s}e^{-r^2/4}$).
Let $\chi_k(r)$ be a cut-off with $\chi_k=1$ for $r\le k$, $\chi_k=0$ for $r\ge k+1$, $|\nabla\chi_k|\le2$, and $\theta_j(t)=\min\{1,(jt-1)_+\}$, so that $\theta_j\uparrow1$ on $(0,\infty)$ and
$0\le\theta_j'\le j\mathbf 1_{(1/j,2/j)}$. The function $\varphi=h\chi_k^2\theta_j(\eta_1)^2$ is Lipschitz with compact support in $\Rp$, hence an admissible test function for the classical equation (Lemma~\ref{lem:test}(a)):
\[
\int\chi_k^2\theta_j^2|\nabla h|^2K+\int\eta_1^{2\alpha}h^{q+1}\chi_k^2\theta_j^2K
=\lambda\int h^2\chi_k^2\theta_j^2K-2\int h\chi_k\theta_j^2\,\nabla h\!\cdot\!\nabla\chi_k\,K-2\int h\chi_k^2\theta_j\theta_j'\,\partial_1h\,K .
\]
Apply $2ab\le\frac14a^2+4b^2$ to each cross term (with $a=\chi_k\theta_j|\nabla h|$, and $b=\theta_jh|\nabla\chi_k|$, resp.\ $b=\chi_k\theta_j'h$):
\[
\tfrac12\int\chi_k^2\theta_j^2|\nabla h|^2K+\int\eta_1^{2\alpha}h^{q+1}\chi_k^2\theta_j^2K\ \le\ \lambda\int h^2K+16\int_{\{k<r<k+1\}}h^2K+4\int\chi_k^2\theta_j'^2h^2K .
\]
On the support of $\theta_j'$ we have $\eta_1<2/j$, so $h\le2C_{k+1}/j$ by Step 1 (with $R=k+1$), and $\theta_j'^2h^2\le4C_{k+1}^2$; hence the last integral is at most
$4C_{k+1}^2\max_{B_{k+1}}K\cdot|\{\frac1j<\eta_1<\frac2j\}\cap B_{k+1}|\to0$ as $j\to\infty$. Letting $j\to\infty$ (monotone convergence) and then $k\to\infty$,
\[
\tfrac12\int_{\Rp}|\nabla h|^2K+\int_{\Rp}\eta_1^{2\alpha}h^{q+1}K\le(\lambda+16)\int_Ah^2K<\infty .
\]

\emph{$h\in W^{1,2}_{0,K}(\Rp)$.} For $\varepsilon>0$ the function $(h-\varepsilon)_+\chi_k$ has compact support in $\Rp$, because $h$ is continuous on $\overline{\Rp}$ and vanishes on $\partial\Rp$; so it lies in $W^{1,2}_{0,K}(\Rp)$.
As $\varepsilon\to0$ and then $k\to\infty$ it converges to $h$ in $W^{1,2}_K$ (dominated convergence, using $h,\nabla h\in L^2_K$ and $h>0$). Hence $h\in W^{1,2}_{0,K}(\Rp)$, and so $h\in\X$.

\emph{The identity \eqref{eq:weak}.} Since $h\in C^2(\Rp)\cap\X$ solves the equation classically, Lemma~\ref{lem:test}(c) gives \eqref{eq:weak} for all $\varphi\in\X$; so $h\in\Hh$.
\end{proof}

\begin{corollary}\label{cor:S}
If $q\ge q_c$ then $\mathcal S=\emptyset$. If $q<q_c$ then $\mathcal S=\{h\}$, where $h$ is the function of Proposition~\ref{prop:h}.
\end{corollary}
\begin{proof}
By Proposition~\ref{prop:S}, $\mathcal S\subset\Hh\setminus\{0\}$. For $q\ge q_c$, $\Hh=\{0\}$ by Proposition~\ref{prop:nonex}. For $q<q_c$, Proposition~\ref{prop:uniq} gives $\mathcal S\subset\{h\}$, and $h\in\mathcal S$ by Proposition~\ref{prop:h} (the bound $h\le CU$ gives \eqref{eq:decay}).
\end{proof}

%---------------------------------------------------------------------
\section{Proof of Theorem~\ref{thm:main}}
(i) and the equalities in (ii) are Proposition~\ref{prop:nonex}, Proposition~\ref{prop:uniq} and Corollary~\ref{cor:S}. Existence and regularity in (ii) are Proposition~\ref{prop:h}, and
\eqref{eq:sharp} follows from Lemma~\ref{lem:U} with Proposition~\ref{prop:h} (upper bound; recall $\sigma_1r^s=\eta_1r^{s-1}$) and Proposition~\ref{prop:lower} (lower bound).

(iii) By Section~2.1, $u_s=t^{-\lambda}h(x/\sqrt t)$ solves \eqref{eq:pde} in $E$, and it is positive and continuous up to $\{x_1=0\}$ with value $0$. From
\eqref{eq:sharp} and $h\le C_1\min(\eta_1,1)$, for $x\ne0$,
\[
u_s(x,t)\le C|x|^{s}\,t^{\frac N2-2\lambda}e^{-|x|^2/(4t)}\to0\quad(t\to0)
\]
locally uniformly in $x\ne0$. If $u$ is any positive self-similar solution in $C(\overline E\setminus\{(0,0)\})$, its profile lies in $\mathcal S$ (Section~\ref{sec:cont}), so $u=u_s$ for $q<q_c$ and no such $u$ exists for $q\ge q_c$ by Corollary~\ref{cor:S}.
Finally $\int_Hx_1u_s\dd x=t^{\frac{N+1}2-\lambda}\int_{\Rp}\eta_1h\dd\eta$ by the change of variables $x=\sqrt t\,\eta$, and
$\int_{\Rp}\eta_1h<\infty$ by the bounds on $h$. Since $\lambda>\frac{N+1}2$, the exponent is negative, so the mass of $x_1u_s$ tends to $+\infty$ as $t\to0$. \qed

%---------------------------------------------------------------------
\section{The asymptotic profile}\label{sec:asym}

Throughout this section $1<q<q_c$ and $h$ is the solution of Theorem~\ref{thm:main}. Put $H(r,\sigma)=h(r\sigma)\,r^{-s}e^{r^2/4}$; by \eqref{eq:sharp}, $c\sigma_1\le H\le C\sigma_1$ for $r\ge R_1$.

\paragraph{Angular modes.} A Dirichlet eigenfunction of $-\Delta_{\Sph}$ on $\Sph^{N-1}_+$ extends, by odd reflection in $\sigma_1$, to a spherical harmonic. Hence the eigenvalues are $\mu_\ell=\ell(\ell+N-2)$, $\ell\ge1$ \cite{BGM}, and the eigenfunctions are the restrictions of spherical harmonics of degree $\ell$ that are odd in $\sigma_1$. Let $(\phi_k)$ be an orthonormal basis of $L^2(\Sph^{N-1}_+)$ of such eigenfunctions, of degrees $\ell_k$, with $\phi_1=\sigma_1/\|\sigma_1\|$, and let $Y_k(\eta)=|\eta|^{\ell_k}\phi_k(\eta/|\eta|)$, a harmonic homogeneous polynomial. Put
\[
A_\ell=\frac{N+\ell}2-\lambda,\qquad b_\ell=\ell+\frac N2,\qquad P(r)=r^{N-1}e^{r^2/4}.
\]
$M(A,b,z)$ and $U(A,b,z)$ denote Kummer's functions \cite[Ch.~13]{DLMF}.

\begin{lemma}[Kummer reduction]\label{lem:kummer}
Let $\ell\ge1$, $\mu=\mu_\ell$, $A=A_\ell$, $b=b_\ell$ and $\mathcal La=a''+(\frac{N-1}r+\frac r2)a'+(\lambda-\frac\mu{r^2})a=P^{-1}(Pa')'+(\lambda-\frac\mu{r^2})a$.
\begin{enumerate}[label=\textup{(\roman*)},leftmargin=*]
\item For $F\in C^2(0,\infty)$: $\mathcal L\big[r^\ell e^{-r^2/4}F(r^2/4)\big]=r^\ell e^{-r^2/4}\big[zF''+(b-z)F'-AF\big]_{z=r^2/4}$.
\item $\upsilon=r^\ell e^{-r^2/4}U(A,b,r^2/4)$ and $\varrho=r^\ell e^{-r^2/4}M(A,b,r^2/4)$ solve $\mathcal L=0$. As $r\to\infty$, $\upsilon=4^Ar^se^{-r^2/4}(1+o(1))$, and $\upsilon>0$ on some $[r_0,\infty)$. Near $0$, $\varrho\sim r^\ell$.
\item $\varpi=\upsilon\int_{r_0}^r\frac{\dd t}{P\upsilon^2}$ solves $\mathcal L=0$, $P(\upsilon\varpi'-\upsilon'\varpi)=1$, and $\varpi=2\cdot4^{-A}r^{-2\lambda}(1+o(1))$.
\item $P(\upsilon\varrho'-\upsilon'\varrho)\equiv w_\ell:=\dfrac{\Gamma(b)4^b}{2\Gamma(A)}$; in particular $w_\ell=0$ if $A\in\{0,-1,-2,\dots\}$.
\end{enumerate}
\end{lemma}
\begin{proof}
(i) is a direct computation. (ii) $U$ and $M$ solve Kummer's equation, and $U(A,b,z)=z^{-A}(1+o(1))$ \cite[13.7.3]{DLMF}. Since $r^\ell(r^2/4)^{-A}=4^Ar^{s}$, this gives the asymptotics of $\upsilon$, and $M(A,b,0)=1$ gives $\varrho\sim r^\ell$. (iii) This is reduction of order. Since $1/(P\upsilon^2)=4^{-2A}r^{1-N-2s}e^{r^2/4}(1+o(1))$, l'H\^opital's rule gives $\int^r\frac{\dd t}{P\upsilon^2}=2\cdot4^{-2A}r^{-N-2s}e^{r^2/4}(1+o(1))$, and $N+s=2\lambda$. (iv) Use the Wronskian $\mathcal W\{M,U\}(z)=-\Gamma(b)z^{-b}e^z/\Gamma(A)$ \cite[13.2.34]{DLMF} (there stated for $\mathbf M=M/\Gamma(b)$) and the chain rule with $z=r^2/4$.
\end{proof}

\begin{proposition}[mode-wise limits]\label{prop:modes}
For every $k$, $a_k(r)=\int_{\Sph^{N-1}_+}h(r\sigma)\phi_k(\sigma)\dd\sigma$ satisfies $\mathcal L a_k=N_k$ on $(0,\infty)$ with $\ell=\ell_k$ and $N_k(r)=\int r^{2\alpha}\sigma_1^{2\alpha}h(r\sigma)^q\phi_k\dd\sigma$. The limit
\[
c_k=\lim_{r\to\infty}\int H(r,\sigma)\phi_k(\sigma)\dd\sigma=\lim_{r\to\infty}4^{A_{\ell_k}}\frac{a_k(r)}{\upsilon(r)}
\]
exists. Consequently $H(r,\cdot)\rightharpoonup\Phi:=\sum_kc_k\phi_k$ weakly in $L^2(\Sph^{N-1}_+)$, $c\,\sigma_1\le\Phi\le C\,\sigma_1$ a.e., and $c_1>0$.
\end{proposition}
\begin{proof}
\emph{The mode equation.} Since $\nabla h\in L^2_K$ and $h\in C(\overline{\Rp})$ vanishes on $\partial\Rp$, for a.e.\ $r$ we have $h(r\cdot)\in W^{1,2}_0(\Sph^{N-1}_+)$, so $\int\nabla_{\Sph}h\cdot\nabla_{\Sph}\phi_k=\mu_k\int h\phi_k$. Testing \eqref{eq:weak} with $\chi(r)\phi_k(\sigma)\in\X$, $\chi\in C_c^\infty(0,\infty)$, gives $\mathcal La_k=N_k$ in $\mathcal D'(0,\infty)$. Here $|\phi_k|\le C_k\sigma_1$, and for $r\ge R_0$, $h\le C\sigma_1r^se^{-r^2/4}$; with $\int\sigma_1^{2\alpha+q+1}\dd\sigma<\infty$ (as $2\alpha+q+1>0$) this gives $|N_k|\le C_kr^{2\alpha+sq}e^{-qr^2/4}$, and $N_k$ is continuous, so $a_k\in C^2$.

\emph{Variation of constants.} On $[r_0,\infty)$ write $a_k=\mathfrak a\upsilon+\mathfrak b\varpi$ with $\mathfrak a'=-P\varpi N_k$ and $\mathfrak b'=P\upsilon N_k$ (Lemma~\ref{lem:kummer}(iii)). Both right-hand sides are $O(r^\kappa e^{-(q-1)r^2/4})$, so $\mathfrak a\to\mathfrak a_\infty$ and $\mathfrak b\to\mathfrak b_\infty$. Since $|a_k|,|\mathfrak a\upsilon|\le Cr^se^{-r^2/4}$, also $|\mathfrak b\varpi|\le Cr^se^{-r^2/4}$, hence $\mathfrak b\to0$. Then $\mathfrak b(r)=-\int_r^\infty P\upsilon N_k=O(r^{\kappa}e^{-qr^2/4})$, so $\mathfrak b\varpi/\upsilon\to0$, and $a_k/\upsilon\to\mathfrak a_\infty$. Since $\upsilon r^{-s}e^{r^2/4}\to4^{A}$, we get $c_k=4^{A}\mathfrak a_\infty$.

\emph{Limits.} The family $H(r,\cdot)$ is bounded in $L^2$ and all its coefficients converge; this gives weak convergence. Testing against nonnegative functions preserves the bounds $c\sigma_1\le\cdot\le C\sigma_1$, and $c_1=\lim\int H\phi_1\ge c\int\sigma_1\phi_1>0$.
\end{proof}

\begin{proposition}[structure of $\Phi$]\label{prop:structure}
\begin{enumerate}[label=\textup{(\roman*)},leftmargin=*]
\item \emph{Axial symmetry.} $h(\eta)$ depends only on $(r,\sigma_1)$ (for $\alpha=0$ this is observed in \cite[Remark after Lemma~2.5]{MV01}). Hence $c_k=0$ unless $\phi_k$ is zonal, and
\[
\Phi(\sigma)=\sum_{\ell\ \mathrm{odd}}c_\ell Z_\ell(\sigma_1),
\]
where $Z_\ell$ is the $L^2$-normalised zonal harmonic of degree $\ell$: a multiple of the Gegenbauer polynomial $C_\ell^{(N-2)/2}$ if $N\ge3$, and of the Chebyshev polynomial $T_\ell$ if $N=2$.
\item \emph{Coefficient formula.} Let $\ell$ be odd and $Y_\ell(\eta)=|\eta|^\ell Z_\ell(\sigma_1)$. If $A_\ell\notin\{0,-1,-2,\dots\}$,
\[
c_\ell=-\frac{2\,\Gamma(A_\ell)}{\Gamma(b_\ell)\,4^{\,b_\ell-A_\ell}}\int_{\Rp}\eta_1^{2\alpha}h(\eta)^q\,Y_\ell(\eta)\,M\Big(A_\ell,b_\ell,\frac{|\eta|^2}4\Big)\dd\eta .
\]
If $A_\ell=-j\in\{0,-1,\dots\}$, the same integral vanishes; there $M(-j,b_\ell,\cdot)$ is a Laguerre polynomial.
\item \emph{Separable case.} If $2\alpha+q=1$, then $h=\sigma_1f(r)$, and $H(r,\cdot)\to c_1\phi_1$ uniformly, so $\Phi$ is a multiple of $\sigma_1$. If $2\alpha+q\ne1$, $h$ is not of the form $\sigma_1f(r)$.
\end{enumerate}
\end{proposition}
\begin{proof}
(i) For every orthogonal map $O$ of $\R^{N-1}$, $(\eta_1,\eta')\mapsto(\eta_1,O\eta')$ preserves $\Rp$, $K$ and $\eta_1$, so $h\circ(\mathrm{id}\times O)\in\Hh\setminus\{0\}$. By Proposition~\ref{prop:uniq} it equals $h$. A function of $\sigma_1$ alone is orthogonal to every non-zonal spherical harmonic. Zonal harmonics of degree $\ell$ have parity $(-1)^\ell$ in $\sigma_1$, so only odd $\ell$ occur.

(ii) Let $\varrho$ be as in Lemma~\ref{lem:kummer}. Since $\mathcal L\varrho=0$ and $\mathcal La_k=N_k$, we have $(P(a_k'\varrho-a_k\varrho'))'=P\varrho N_k$; integrate over $(\epsilon,R)$.

\emph{At $R\to\infty$.} With $a_k=\mathfrak a\upsilon+\mathfrak b\varpi$, $P(a_k'\varrho-a_k\varrho')=-\mathfrak aw_\ell+\mathfrak b\,P(\varpi'\varrho-\varpi\varrho')$. The last Wronskian is constant and $\mathfrak b\to0$, so the limit is $-\mathfrak a_\infty w_\ell$.

\emph{At $\epsilon\to0$.} We have $|a_k|\le Cr$ (Lemma~\ref{lem:bdry}), $\varrho\sim r^\ell$, and $\int_0^1|a_k'|^2r^{N-1}<\infty$ (because $\nabla h\in L^2_K$). Hence the term tends to $0$ along a sequence $\epsilon_j\to0$.

$P\varrho N_k$ is integrable on $(0,\infty)$ (near $0$ because $2\alpha+q+1>0$, near $\infty$ by the Gaussian bound), and
\[
\int_0^\infty P\varrho N_k\dd r=\int_{\Rp}\eta_1^{2\alpha}h^qY_\ell M\dd\eta,\qquad\text{so}\qquad -\mathfrak a_\infty w_\ell=\int_{\Rp}\eta_1^{2\alpha}h^qY_\ell M,
\] which with $c_\ell=4^{A_\ell}\mathfrak a_\infty$ and Lemma~\ref{lem:kummer}(iv) gives both claims.

(iii) Let $2\alpha+q=1$ and let $\Pi$ be the projection onto the first angular mode, $\Pi\varphi=\sigma_1\,g(r)$ with $g=\int\varphi\sigma_1\dd\sigma/\int\sigma_1^2\dd\sigma$.
\begin{itemize}[leftmargin=*]
\item $\Pi$ maps $X_n$ into itself. It preserves $W_0^{1,2}(D_n)$, and $\int\eta_1^{2\alpha}|\Pi\varphi|^{q+1}\le C\int\eta_1^{2\alpha}|\varphi|^{q+1}$ by H\"older on the sphere, since $(q+1-2\alpha)/q=2>-1$.
\item Minimise $J$ over $\Pi X_n$. As in Proposition~\ref{prop:hn} this gives $\tilde h_n=\sigma_1f_n\ge0$, $\tilde h_n\ne0$, with $J'(\tilde h_n)[\Pi\varphi]=0$.
\item Since $\eta_1^{2\alpha}\tilde h_n^q=\sigma_1r^{2\alpha}f_n^q$, each term of $J'(\tilde h_n)[\varphi-\Pi\varphi]$ is a radial integral of $\int(\varphi-\Pi\varphi)\sigma_1\dd\sigma=0$ or of $\int\nabla_{\Sph}\sigma_1\cdot\nabla_{\Sph}(\varphi-\Pi\varphi)\dd\sigma=(N-1)\int\sigma_1(\varphi-\Pi\varphi)\dd\sigma=0$.
\item Hence $\tilde h_n$ is a nonnegative nontrivial weak solution in $X_n$, positive by the strong maximum principle. Lemma~\ref{lem:comp}, applied in both directions, gives $\tilde h_n=h_n$. Letting $n\to\infty$, $h=\sigma_1f$.
\end{itemize}
Conversely, if $h=\sigma_1f$, then $\eta_1^{2\alpha}h^q=r^{2\alpha}\sigma_1^{2\alpha+q}f^q$ must be a multiple of $\sigma_1$ on each sphere, which forces $2\alpha+q=1$.
\end{proof}

\begin{lemma}[angular integrals]\label{lem:sign}
For odd $\ell$ and $\gamma>-1$, let $S_\ell(\gamma)=\int_{\Sph^{N-1}_+}\sigma_1^\gamma Z_\ell(\sigma_1)\dd\sigma$, with $Z_\ell$ normalised positive at the pole. Then
\[
S_\ell(\gamma)=\frac{p_{\ell,N}(\gamma)}{\Gamma\big(\frac{\gamma+2-\ell}2\big)},\qquad p_{\ell,N}(\gamma)>0 .
\]
Explicitly, with $\nu=\frac{N-2}2$ and $\mu=\gamma+1$, for $N\ge3$
\[
\int_0^1t^\gamma C_\ell^\nu(t)(1-t^2)^{\nu-\frac12}\dd t=\frac{\pi\,2^{1-2\nu-\mu}\,\Gamma(\ell+2\nu)\,\Gamma(\mu)}{\ell!\,\Gamma(\nu)\,\Gamma(\frac{1+\mu+\ell}2+\nu)\,\Gamma(\frac{1+\mu-\ell}2)},
\]
and for $N=2$, $\int_{-\pi/2}^{\pi/2}\cos^\gamma\!\varphi\cos\ell\varphi\dd\varphi=\frac{\pi\Gamma(\gamma+1)}{2^\gamma\Gamma(1+\frac{\gamma+\ell}2)\Gamma(\frac{\gamma+2-\ell}2)}$. In particular $S_3(\gamma)$ has the sign of $\gamma-1$, and $S_\ell(\gamma)$ vanishes exactly when $\gamma+2-\ell\in\{0,-2,-4,\dots\}$.
\end{lemma}
\begin{proof}
For $N\ge3$ this is \cite[7.311.2]{GR} with $n=\ell$ and $n+2\rho=\gamma$, rewritten with the duplication formula for $\Gamma$. There it is stated for $\rho>-\frac12$, i.e.\ $\gamma>\ell-1$; both sides are analytic in $\gamma$ on $\{\operatorname{Re}\gamma>-1\}$, so the identity holds for all $\gamma>-1$. For $N=2$ it is \cite[3.631.9]{GR} with $\nu=\gamma+1$, $a=\ell$. For $N\ge3$, $\int_{\Sph^{N-1}_+}f(\sigma_1)\dd\sigma=|\Sph^{N-2}|\int_0^1f(t)(1-t^2)^{\nu-\frac12}\dd t$, and $Z_\ell$ is a positive multiple of $C_\ell^\nu$. For $\ell=3$: $(\gamma-1)/2\in(-1,0)$ iff $\gamma\in(-1,1)$, where $\Gamma<0$.
\end{proof}

\begin{theorem}[near-critical shape]\label{thm:nearcrit}
Fix $N\ge2$, $\alpha>-1$, and put $\gamma_c=2\alpha+q_c$. For every odd $\ell\ge3$, as $q\uparrow q_c$ (so $\lambda\downarrow\lambda_1$),
\[
\frac{c_\ell}{c_1}=-(\lambda-\lambda_1)\,\kappa_\ell+o(\lambda-\lambda_1),\qquad
\kappa_\ell=\frac{\Gamma(b_1)\Gamma(A_\ell^0)\,4^{\,b_1-b_\ell+A_\ell^0}}{\Gamma(b_\ell)}\cdot\frac{\mathcal R_\ell}{\mathcal R_1}\cdot\frac{S_\ell(\gamma_c)}{S_1(\gamma_c)},
\]
where $A_\ell^0=\frac{\ell-1}2$ and $\mathcal R_\ell=\int_0^\infty r^{\gamma_c+N+\ell-1}e^{-q_cr^2/4}M(A_\ell^0,b_\ell,\frac{r^2}4)\dd r>0$. Consequently, for each odd $\ell\ge3$ with $\gamma_c+2-\ell\notin\{0,-2,-4,\dots\}$ there is $q_\ell\in(1,q_c)$ such that, for every $q\in(q_\ell,q_c)$,
\[
\operatorname{sign}c_\ell=-\operatorname{sign}\frac1{\Gamma\big(\frac{\gamma_c+2-\ell}2\big)} .
\]
The threshold $q_\ell$ may depend on $\ell$; no statement uniform in $\ell$ is made. In particular:
\begin{enumerate}[label=\textup{(\alph*)},leftmargin=*]
\item if $\gamma_c<3$ (i.e.\ $\alpha<\frac N{N+2}$) and $\gamma_c\ne1$ (i.e.\ $\alpha\ne-\frac1{N+2}$), then for every odd $\ell\ge3$, $\operatorname{sign}c_\ell=(-1)^{\frac{\ell-1}2}\operatorname{sign}(\gamma_c-1)$ for $q\in(q_\ell,q_c)$. Hence, for any finite set of odd degrees, the signs alternate simultaneously when $q$ is close enough to $q_c$;
\item if $\gamma_c>3$ and $\gamma_c$ is not an odd integer, then $c_\ell<0$ for all odd $\ell$ with $3\le\ell<\gamma_c+2$, simultaneously, when $q$ is close enough to $q_c$; in particular $c_3$ and $c_5$ are both negative, so the signs do not alternate from the first two modes. For each odd $\ell>\gamma_c+2$, $c_\ell$ and $c_{\ell+2}$ have opposite signs for $q$ close enough to $q_c$ (depending on $\ell$).
\end{enumerate}
\end{theorem}
\begin{proof}
Let $\bar h=h/\|h\|_{L^2_K}$.

\emph{Step 1: $\bar h\to\hat\psi_1$.} The energy identity gives $\int|\nabla\bar h|^2K\le\lambda$. Write $\bar h=a\hat\psi_1+g$ with $g\perp\psi_1$ in $L^2_K$; then $g\in W^{1,2}_{0,K}$ and, by Lemma~\ref{lem:spec}(ii), $\lambda\ge\int|\nabla\bar h|^2K=\lambda_1a^2+\int|\nabla g|^2K\ge\lambda_1+(\lambda_2-\lambda_1)\|g\|^2_{L^2_K}$. So $g$ has squared norm at most $(\lambda-\lambda_1)/(\lambda_2-\lambda_1)\to0$. Since $\bar h\ge0$, $\bar h\to\hat\psi_1=\psi_1/\|\psi_1\|_{L^2_K}$ in $L^2_K$.

\emph{Step 2: uniform bounds.} Restrict to $\lambda\in(\lambda_1,\lambda_1+1]$. The function $\bar h$ satisfies $(L_K-\lambda)\bar h+\epsilon\,\eta_1^{2\alpha}\bar h^q=0$ with $\epsilon=\|h\|^{q-1}_{L^2_K}>0$, and $\|\bar h\|_{L^2_K}=1$. We check that each bound used below has constants independent of $\lambda$ in this range.
\begin{itemize}[leftmargin=*]
\item \emph{Sup bound.} The proof of Lemma~\ref{lem:reg} uses only $\bar h\ge0$ and the linear part of the equation (the absorption term is dropped because it is nonnegative). It gives $\sup_{B_1(y)\cap\Rp}\bar h\le C(N,\lambda)e^{-(|y|-1)_+^2/8}$, and $C(N,\lambda)$ comes from Moser iteration for $-\Delta\hat v\le\lambda\hat v$, so it is bounded for $\lambda\le\lambda_1+1$.
\item \emph{Boundary bound.} The barrier $\Theta$ of Lemma~\ref{lem:bdry} satisfies $(L_K-\lambda)\Theta\ge0$ by a linear computation, so it is a supersolution for every $\epsilon\ge0$. Its constants are $m=\sup_{|\eta|\le R+2}\bar h$, bounded by the previous item, and $\rho_*(N,R,\lambda)$, which is bounded below for $\lambda\le\lambda_1+1$. Hence $\bar h\le C_R\eta_1$ on $\{|\eta|\le R\}$ with $C_R$ independent of $\lambda$.
\item \emph{Gaussian bound.} In Lemmas~\ref{lem:U} and~\ref{lem:Z} the inequalities $(L_K-\lambda)U\ge0$ and $(L_K-\lambda)Z>0$ are linear, with $R_0,R_0'$ depending continuously on $\lambda$. Step~2 of Proposition~\ref{prop:S} then gives $\bar h\le C\,U_\lambda$ in $\{r\ge R_0\}$, with $C=2C_{R_0+1}\sup_{R_0\le r\le R_0+1}r^{1-s}e^{r^2/4}$, again independent of $\lambda$.
\end{itemize}
Together: $\bar h\le C\min(\eta_1,1)$ in $\{r\le R_0\}$ and $\bar h\le C\sigma_1r^se^{-r^2/4}$ in $\{r\ge R_0\}$, uniformly for $\lambda\in(\lambda_1,\lambda_1+1]$. Moreover $|M(A,b,z)|\le e^z$ for $|A|\le1$, $b\ge1$, and also for $0<A\le b$ (then $(A)_n\le(b)_n$).

\emph{Step 3: the limit.} Write $c_\ell=-\frac{2\Gamma(A_\ell)}{\Gamma(b_\ell)4^{b_\ell-A_\ell}}\|h\|^q\mathcal I_\ell$ with $\mathcal I_\ell=\int\eta_1^{2\alpha}\bar h^qY_\ell M(A_\ell,b_\ell,\cdot)$ (Proposition~\ref{prop:structure}(ii)). Here $A_1=\lambda_1-\lambda\in(-1,0)$ and $A_\ell\to A_\ell^0>0$ for $\ell\ge3$. By Step 2 and dominated convergence, $\mathcal I_\ell\to\int\eta_1^{2\alpha}\hat\psi_1^{q_c}Y_\ell M(A_\ell^0,b_\ell,\cdot)$, with $A_1^0=0$ (so $M\equiv1$). Since $\eta_1^{2\alpha}\psi_1^{q_c}=r^{\gamma_c}\sigma_1^{\gamma_c}e^{-q_cr^2/4}$, these integrals factor into $\mathcal R_\ell\,S_\ell(\gamma_c)$, up to a common positive constant. Dividing the formulas for $c_\ell$ and $c_1$, and using $\Gamma(A_1)=\Gamma(1+A_1)/A_1=-\Gamma(1+A_1)/(\lambda-\lambda_1)$, we get $c_\ell/c_1=-(\lambda-\lambda_1)(\kappa_\ell+o(1))$; this is the stated expansion, and it determines the sign of $c_\ell$ when $\kappa_\ell\ne0$. Also and $\mathcal R_\ell>0$ because $M(A,b,\cdot)>0$ for $A,b>0$. The sign law follows from Lemma~\ref{lem:sign}.

\emph{The two cases.} For every odd $\ell\ge3$ with $\kappa_\ell\ne0$, the expansion gives
\[
\lim_{q\uparrow q_c}\frac{c_\ell}{c_1(\lambda-\lambda_1)}=-\kappa_\ell .
\]
Since $c_1>0$ and $\lambda>\lambda_1$, $\operatorname{sign}c_\ell=-\operatorname{sign}\kappa_\ell$ for $q$ in a left neighbourhood $(q_\ell,q_c)$ of $q_c$, which may depend on $\ell$. For a finite set of degrees, a common neighbourhood is the intersection of the corresponding ones. All factors of $\kappa_\ell$ other than $S_\ell(\gamma_c)$ are positive ($A_\ell^0=\frac{\ell-1}2>0$, $\mathcal R_\ell>0$, and $S_1(\gamma_c)>0$ since $\gamma_c>-1$), so Lemma~\ref{lem:sign} gives the sign; $\kappa_\ell\ne0$ exactly when $x_\ell=\frac{\gamma_c+2-\ell}2$ is not a nonpositive integer, which holds for all odd $\ell\ge3$ when $\gamma_c$ is not an odd integer. For $x\in(-k-1,-k)$, $\operatorname{sign}\Gamma(x)=(-1)^{k+1}$. If $\gamma_c<3$, then $x_\ell<0$ for $\ell\ge5$, and $x_\ell$ decreases by $1$ when $\ell$ increases by $2$, so the signs of $1/\Gamma(x_\ell)$ alternate; this gives (a), the case $\ell=3$ being Lemma~\ref{lem:sign}. If $\gamma_c>3$, then $x_\ell>0$ for the finitely many odd $\ell<\gamma_c+2$, so these coefficients are simultaneously negative; for odd $\ell>\gamma_c+2$ consecutive factors $1/\Gamma(x_\ell)$ have opposite signs. This gives (b).
\end{proof}

\paragraph{Exceptional values.}
By Lemma~\ref{lem:sign}, $\kappa_\ell=0$ exactly when $\gamma_c$ is an odd integer and $\ell\ge\gamma_c+2$; there Theorem~\ref{thm:nearcrit} gives only $c_\ell=o(\lambda-\lambda_1)c_1$, and a higher-order expansion would be needed. The case $\gamma_c=1$, i.e.\ $\alpha=-\frac1{N+2}$, is the one in which $\psi_1$ itself is the separable profile ($2\alpha+q_c=1$): then $\kappa_\ell=0$ for all $\ell\ge3$. In that case $2\alpha+q<1$ for every $q<q_c$, and Theorem~\ref{thm:sign} below gives $c_3>0$ for all $1<q<q_c$. For $\gamma_c=3$, $c_3<0$ near $q_c$ by Theorem~\ref{thm:nearcrit}, while $\kappa_\ell=0$ for $\ell\ge5$.

\subsection*{Convergence, analyticity and monotonicity of the profile}
By Proposition~\ref{prop:structure}(i) only the odd zonal modes occur. For odd $\ell$ write $c_\ell(r)=\int H(r,\cdot)Z_\ell\dd\sigma$, and put $\kappa_1=2s+N-1$, $\omega(r)=r^{\kappa_1}e^{-r^2/4}$, and
\[
\beta_\ell=\mu_\ell-s(s+N-2)=(\ell-s)(\ell+s+N-2),\qquad B_\ell=b_\ell-A_\ell-1=\tfrac\ell2+\lambda-1 .
\]
Note that $A_\ell=(\ell-s)/2$, so $\ell>s$ iff $A_\ell>0$ iff $\beta_\ell>0$.

\begin{lemma}[modes of $H$]\label{lem:cmode}
\begin{enumerate}[label=\textup{(\roman*)},leftmargin=*]
\item $(\omega c_\ell')'-\omega\beta_\ell r^{-2}c_\ell=\omega f_\ell$ on $(0,\infty)$, where $f_\ell=\int F(r,\cdot)Z_\ell\dd\sigma$ and $F=r^{2\alpha+s(q-1)}\sigma_1^{2\alpha}e^{-(q-1)r^2/4}H^q\ge0$.
\item For $\ell>s$, the recessive solution of the homogeneous equation is $\tilde\upsilon_\ell=r^{-s}e^{r^2/4}\upsilon_\ell=4^{A_\ell}V_\ell(r^2/4)$, where
\[
V_\ell(z)=z^{A_\ell}U(A_\ell,b_\ell,z)=\mathbb E\big[(1+X/z)^{B_\ell}\big],\qquad X\sim\mathrm{Gamma}(A_\ell,1).
\]
Hence $V_\ell$ is decreasing, $V_\ell\to1$ as $z\to\infty$, and $V_\ell(z)\ge(1+A_\ell/z)^{B_\ell}$.
\item There are $C,\kappa$ with $\|F(t,\cdot)\|_{H^{-1}(\Sph^{N-1}_+)}\le m(t):=Ct^{\kappa}e^{-(q-1)t^2/4}$ for $t\ge R_1$. Hence $|f_\ell(t)|\le\mu_\ell^{1/2}m(t)$.
\end{enumerate}
\end{lemma}
\begin{proof}
(i) Substitute $a_\ell=r^se^{-r^2/4}c_\ell$ in the mode equation of Proposition~\ref{prop:modes}.

(ii) For $A>0$, $U(A,b,z)=\Gamma(A)^{-1}\int_0^\infty e^{-zt}t^{A-1}(1+t)^{b-A-1}\dd t$ \cite[13.4.4]{DLMF}; substitute $t=x/z$. Since $B_\ell\ge1$, Jensen's inequality gives $V_\ell(z)\ge(1+\mathbb EX/z)^{B_\ell}=(1+A_\ell/z)^{B_\ell}$.

(iii) For $t\ge R_1$, $F\le Ct^{2\alpha+sq-s}e^{-(q-1)t^2/4}\sigma_1^{2\alpha+q}$. For $v\in H^1_0(\Sph^{N-1}_+)$, Hardy's inequality $\int v^2\sigma_1^{-2}\le C\int|\nabla_{\Sph}v|^2$ and Cauchy--Schwarz give $\int\sigma_1^{2\alpha+q}|v|\le\big(\int\sigma_1^{2(2\alpha+q+1)}\big)^{1/2}C\|\nabla_{\Sph}v\|_{L^2}$, and $2(2\alpha+q+1)>-1$. Finally $\|Z_\ell\|_{H^1_0}=\mu_\ell^{1/2}$.
\end{proof}

\begin{theorem}[convergence and analyticity]\label{thm:conv}
\begin{enumerate}[label=\textup{(\roman*)},leftmargin=*]
\item $H(r,\cdot)\to\Phi$ in $L^p(\Sph^{N-1}_+)$ for every $p<\infty$.
\item There are $C,\theta>0$ with $|c_\ell|\le C(1+\ell)e^{-\theta\ell}$; one may take $\theta=\min\{\frac12\ln2,\frac{q-1}8\}$.
\item $\Phi=\mathcal F(\sigma_1)$, where $\mathcal F$ is an odd function, holomorphic in a complex neighbourhood of $[-1,1]$. In particular $\Phi$ is the restriction of a real-analytic function on $\Sph^{N-1}$, and $\Phi\in C^\infty(\overline{\Sph^{N-1}_+})$.
\end{enumerate}
\end{theorem}
\begin{proof}
Fix $R\ge\max\{R_1,2\sqrt{\kappa_1},\sqrt{32}\}$ and $\ell>s$; write $\tilde\upsilon=\tilde\upsilon_\ell$ and $\mathcal L c=(\omega c')'-\omega\beta_\ell r^{-2}c$.

\emph{Step 1 (a growing solution).} $D=\tilde\upsilon\int_R^r\frac{\dd\tau}{\omega\tilde\upsilon^2}$ satisfies $\mathcal LD=0$, $D(R)=0$, $D>0$ on $(R,\infty)$, $D\to\infty$, and $\omega(\tilde\upsilon D'-\tilde\upsilon'D)=1$.

\emph{Step 2 (maximum principle).} If $v$ is bounded, $\mathcal Lv\le0$ on $(R,\infty)$ and $v(R)\ge0$, then $v\ge0$. Indeed, $v+\delta D$ cannot have a negative interior minimum (there $\mathcal L(v+\delta D)\ge-\omega\beta_\ell r^{-2}(v+\delta D)>0$), it is $\ge0$ at $R$, and it tends to $+\infty$; let $\delta\to0$.

\emph{Step 3 (tails).} Write $c_\ell=c_\ell(R)\tilde\upsilon/\tilde\upsilon(R)+S_\ell$, so $S_\ell$ is bounded, $S_\ell(R)=0$ and $\mathcal LS_\ell=\omega f_\ell$. The constant $M_\ell=\sup_{t\ge R}t^2|f_\ell(t)|/\beta_\ell$ satisfies $\mathcal L(M_\ell\mp S_\ell)\le-\omega|f_\ell|\mp\omega f_\ell\le0$, so by Step 2 $|S_\ell|\le M_\ell\le\mu_\ell^{1/2}\beta_\ell^{-1}E_R$, with $E_R=\sup_{t\ge R}t^2m(t)$. As $0<\tilde\upsilon(r)/\tilde\upsilon(R)\le1$,
\[
\sum_{\ell>L}c_\ell(r)^2\le2\sum_{\ell>L}c_\ell(R)^2+2E_R^2\sum_{\ell>L}\frac{\mu_\ell}{\beta_\ell^2}\qquad(r\ge R),
\]
and both series converge: there is one zonal mode per odd $\ell$, and $\mu_\ell/\beta_\ell^2=O(\ell^{-2})$. So the tails are small uniformly in $r$. Together with $c_\ell(r)\to c_\ell$ for every $\ell$ (Proposition~\ref{prop:modes}), this gives $H(r,\cdot)\to\Phi$ in $L^2$. Since $|H|,|\Phi|\le C\sigma_1$, interpolation gives every $p<\infty$.

\emph{Step 4 (representation).} Variation of constants with $\tilde\upsilon$ and $D$ (Step 1) gives the bounded solution
\[
c_\ell(r)=\tilde\upsilon(r)\Big[\frac{c_\ell(R)}{\tilde\upsilon(R)}-\int_R^rD\omega f_\ell\Big]-D(r)\int_r^\infty\tilde\upsilon\omega f_\ell .
\]
Letting $r\to\infty$, and using $\tilde\upsilon(\infty)=4^{A_\ell}$,
\[
c_\ell=\frac{c_\ell(R)}{V_\ell(R^2/4)}-\int_R^\infty k_\ell f_\ell,\qquad 0\le k_\ell(t)=4^{A_\ell}D(t)\omega(t)\le\frac8t\,V_\ell(t^2/4)^{-1}.
\]
The bound on $k_\ell$ follows because $\tilde\upsilon$ is decreasing and $\omega(t)\int_R^t\omega^{-1}\le8/t$ for $t\ge R$.

\emph{Step 5 (decay).} For $\ell\ge\max\{2s,R^2\}$ we have $A_\ell\ge\ell/4$ and $B_\ell\ge\ell/2$, so $V_\ell(t^2/4)^{-1}\le(1+\ell/t^2)^{-\ell/2}$. The first term is therefore $\le C2^{-\ell/2}$. In the integral, use $(1+\ell/t^2)^{-\ell/2}\le2^{-\ell/2}$ for $t^2\le\ell$, and $m(t)\le Ce^{-(q-1)\ell/8}t^\kappa e^{-(q-1)t^2/8}$ for $t^2>\ell$. This gives (ii).

\emph{Step 6 (analyticity).} $\|Z_\ell\|_{L^\infty}\le C\ell^{(N-2)/2}$. By the Bernstein--Walsh inequality,
\[
|Z_\ell|\le\varrho^\ell\|Z_\ell\|_{L^\infty[-1,1]}\quad\text{on the ellipse }E_\varrho\text{ with foci }\pm1.
\] So $\sum c_\ell Z_\ell$ converges uniformly on $E_\varrho$ for $1<\varrho<e^\theta$.
\end{proof}

\begin{theorem}[monotonicity and the sign of $c_3$]\label{thm:sign}
Let $\gamma=2\alpha+q$ and $u=h/\eta_1$.
\begin{enumerate}[label=\textup{(\roman*)},leftmargin=*]
\item For every $r>0$, $\sigma_1\mapsto u(r\sigma)$ is nonincreasing on $(0,1]$ if $\gamma>1$, nondecreasing if $\gamma<1$, and constant if $\gamma=1$.
\item $\Phi/\sigma_1$ is monotone in $\sigma_1$ in the same sense. If $\gamma\ne1$, it is not constant, and
\[
\operatorname{sign}c_3=-\operatorname{sign}(2\alpha+q-1).
\]
\end{enumerate}
\end{theorem}
\begin{proof}
\emph{Step 1 (lifting to $\R^{N+2}$).} Write $\xi=(y,\eta')\in\R^3\times\R^{N-1}$ and $U(\xi)=u(|y|,\eta')$. Then $U\in H^1_{\rm loc}(\R^{N+2})$: indeed $|\nabla U(\xi)|=|\nabla h-(h/\eta_1)e_1|/\eta_1$ at $(\eta_1,\eta')=(|y|,\eta')$ and $\dd\xi=4\pi\eta_1^2\dd\eta_1\dd\eta'$, so $\int_{B_\rho}|\nabla U|^2\le4\pi\int_{B_\rho\cap\Rp}|\nabla h-(h/\eta_1)e_1|^2<\infty$, since $h/\eta_1$ is locally bounded by Lemmas~\ref{lem:bdry} and~\ref{lem:hopf}. Also $U>0$, $U\le Cr^{s-1}e^{-r^2/4}$ at infinity, and $U$ is a weak solution in $\R^{N+2}$ of
\begin{equation}\label{eq:lift}
-\Delta U-\tfrac12\xi\cdot\nabla U-(\lambda+\tfrac12)U+|y|^{\gamma-1}U^q=0 .
\end{equation}
To see this, average a test function $\Psi\in C_c^\infty(\R^{N+2})$ over $SO(3)$ acting on $y$, giving $\bar\psi(\rho,\eta')$, and apply \eqref{eq:weak} to $\varphi=\eta_1\bar\psi(\eta_1,\eta')\in\X$. The identity $\mathcal M(\eta_1u)=\eta_1\big[-\Delta u-\frac2{\eta_1}\partial_1u-\frac12\eta\cdot\nabla u-(\lambda+\frac12)u+\eta_1^{\gamma-1}u^q\big]$ and $\Delta_{\R^3}F(|y|)=F''+\frac2\rho F'$ turn it into the weak form of \eqref{eq:lift}. The weight $|y|^{\gamma-1}$ is locally integrable because $\gamma>-1$. Finally $U$ is even in $y_1$ and in $\eta_2$.

\emph{Step 2 (rotating planes).} Fix $\beta\in(0,\frac\pi2)$ and let $T$ be the reflection, in the $(y_1,\eta_2)$-plane, across the line at angle $\beta$ (the angle $\vartheta$ goes to $2\beta-\vartheta$); $T$ acts as the identity on the other coordinates. Let $\Omega$ be the set of $\xi$ whose $(y_1,\eta_2)$-component has angle in $(\beta,\beta+\frac\pi2)$. On $\Omega$,
\[
|y(T\xi)|^2-|y(\xi)|^2=\rho_{12}^2\big(\cos^2(\beta-x)-\cos^2(\beta+x)\big)=\rho_{12}^2\sin2\beta\sin2x>0,\qquad x=\vartheta-\beta\in(0,\tfrac\pi2).
\]
$V=U\circ T$ solves \eqref{eq:lift} with weight $|y(T\xi)|^{\gamma-1}$. Hence, if $\gamma>1$, $V$ is a weak subsolution of \eqref{eq:lift} in $\Omega$ (and a supersolution if $\gamma<1$). Moreover $V-U$ has zero trace on $\partial\Omega$: $T$ fixes the ray $\vartheta=\beta$ pointwise, and on the ray $\vartheta=\beta+\frac\pi2$ it acts as $(y_1,\eta_2)\mapsto-(y_1,\eta_2)$, which leaves $U$ invariant. So $V\le U$ in $\Omega$ if $\gamma>1$, and $V\ge U$ if $\gamma<1$, by the comparison principle of Lemma~\ref{lem:comp} transposed to $\R^{N+2}$, which we now justify.

\emph{Comparison in $\R^{N+2}$.} Put $\tilde K(\xi)=e^{|\xi|^2/4}$, $w(\xi)=|y|^{\gamma-1}$, and
$\tilde{\mathcal B}(u;\varphi)=\int\nabla u\cdot\nabla\varphi\,\tilde K-(\lambda+\frac12)\int u\varphi\,\tilde K+\int w\,u^q\varphi\,\tilde K$.
The proof of Lemma~\ref{lem:comp} uses only the Picone identity, the monotonicity of $s\mapsto s^q/(s+\varepsilon)$, and $w>0$ a.e.; it does not use the dimension or the form of the weight. So it applies in $Q=\Omega\subset\R^{N+2}$ to $\tilde{\mathcal B}$, once (H1)--(H3) hold. Take $\gamma>1$, $u_1=V$, $u_2=U$ (for $\gamma<1$ exchange the roles: $u_1=U$, $u_2=V$).
\begin{itemize}[leftmargin=*]
\item \emph{Integrability.} $w\in L^1_{\rm loc}(\R^{N+2})$, since $y\in\R^3$ and $\gamma-1>-2>-3$; this covers the singular case $\gamma<1$. $U$ and $V$ are bounded ($h/\eta_1$ is locally bounded by Lemmas~\ref{lem:bdry} and~\ref{lem:hopf}, and $U\le Cr^{s-1}e^{-r^2/4}$ at infinity), positive, and in $H^1_{\rm loc}(\R^{N+2})$ by Step~1, $T$ being an isometry. Hence $U^2\tilde K$, $V^2\tilde K$ and $w(U^{q+1}+V^{q+1})\tilde K$ are integrable on bounded sets, and $U^2\tilde K\le Cr^{2s-2}e^{-r^2/4}$ is integrable at infinity. This gives (H2).
\item \emph{(H1).} $U$ satisfies $\tilde{\mathcal B}(U;\varphi)=0$ for $\varphi\in C_c^\infty(\R^{N+2})$ by Step~1. For $0\le\varphi\in W^{1,2}_0(\Omega)\cap L^\infty$ with bounded support, take $\varphi_n\in C_c^\infty(\Omega)$, $0\le\varphi_n\le\|\varphi\|_\infty+1$, with $\varphi_n\to\varphi$ in $W^{1,2}$ and a.e. The linear terms converge, and $wU^q\varphi_n\tilde K\to wU^q\varphi\tilde K$ in $L^1$ by dominated convergence, since $wU^q\tilde K$ is integrable on the common bounded support. So $\tilde{\mathcal B}(U;\varphi)=0$. The function $V$ satisfies the same identity with $w\circ T$ in place of $w$. Since $w\circ T\ge w$ in $\Omega$ when $\gamma>1$, $\tilde{\mathcal B}(V;\varphi)=\int(w-w\circ T)V^q\varphi\tilde K\le0$ for $\varphi\ge0$; when $\gamma<1$ the inequality is reversed.
\item \emph{(H3).} $S_\varepsilon=\{u_1>u_2+\varepsilon\}$ is bounded, because $U,V\to0$ at infinity. $\nabla U,\nabla V\in L^2(S_\varepsilon)$ because $U,V\in H^1_{\rm loc}$. Finally, $\Omega$ is the product of an open planar sector with $\R^N$, a Lipschitz domain; $u_1-u_2\in H^1_{\rm loc}(\overline\Omega)$ is continuous on $\partial\Omega$ away from the edge $\{y_1=\eta_2=0\}$, where it vanishes, so its trace on $\partial\Omega$ is zero. Hence $(u_1-u_2-\varepsilon)_+$, which has bounded support, belongs to $W^{1,2}_0(\Omega)$.
\end{itemize}

For points of the $(y_1,\eta_2)$-plane at radius $\rho$, $\sigma_1=|\cos\vartheta|$. Given $0\le\vartheta_1<\vartheta_2<\frac\pi2$, take $\beta=\frac{\vartheta_1+\vartheta_2}2$; this gives $u(\rho\cos\vartheta_1,\cdot)\le u(\rho\cos\vartheta_2,\cdot)$ when $\gamma>1$, which is (i). (Both points avoid $y=0$, where $U$ is continuous.) If $\gamma=1$ the weight is $1$, and (i) is Proposition~\ref{prop:structure}(iii).

\emph{Step 3 (the limit).} By (i), $g_r=H(r,\cdot)/\sigma_1=r^{1-s}e^{r^2/4}u(r\,\cdot)$ is monotone in $\sigma_1$. By Theorem~\ref{thm:conv}, $g_r\to\Phi/\sigma_1$ a.e.\ along a sequence. So $\Phi/\sigma_1$ is monotone.

\emph{Step 4 (not constant).} Let $\gamma>1$ (the case $\gamma<1$ is symmetric) and $W=U-V\ge0$ in $\Omega$. Subtracting the equations of $U$ and $V$,
\[
\mathcal PW:=\big(-\Delta-\tfrac12\xi\cdot\nabla-\lambda-\tfrac12+c\big)W=\big(|y(T\xi)|^{\gamma-1}-|y|^{\gamma-1}\big)V^q\ge0,
\]
where $c=|y|^{\gamma-1}(U^q-V^q)/(U-V)\in\big[0,\,q|y|^{\gamma-1}\max(U,V)^{q-1}\big]$.
The right side is positive in $\Omega$, so $W\not\equiv0$, and the strong maximum principle \cite[Theorem~3.5]{GT} (applied where $y\ne0$, where $U$ and $V$ are $C^2$) gives $W>0$ in $\Omega$.

\smallskip
\emph{(a) Cones and radial solutions.} Fix a direction in $\Omega$ with $y\ne0$, and two open cones $\Omega'\Subset\Omega''\Subset\Omega$ around it whose closures avoid $\{y=0\}$. Let $D'\Subset D''$ be their sections on $\Sph^{N+1}$, and $(\varphi',\mu')$, $(\varphi'',\mu'')$ the first Dirichlet eigenpairs of $-\Delta_{\Sph^{N+1}}$ on $D'$, $D''$; then $\varphi''>0$ on $\overline{D'}$. On $\overline{\Omega''}$, $|y|\ge c_*r$ and $U,V\le Cr^{s-1}e^{-r^2/4}$, so $c\le\bar c(r):=Cr^\kappa e^{-(q-1)r^2/4}$ for some $\kappa$. For $\mu\in\{\mu',\mu''\}$ the radial equation
\[
-\zeta''-\big(\tfrac{N+1}r+\tfrac r2\big)\zeta'+\big(\tfrac{\mu}{r^2}-\lambda-\tfrac12\big)\zeta=0
\]
is the equation of Lemma~\ref{lem:kummer} in dimension $N+2$ with $\lambda+\frac12$ in place of $\lambda$. It has a recessive solution $\zeta_\mu\sim C_\mu r^{s-1}e^{-r^2/4}$, positive on some $[r_1,\infty)$, and a second solution $\varpi_\mu\sim C'_\mu r^{-2\lambda-1}$ (Lemma~\ref{lem:kummer}(ii)--(iii); note $2(\lambda+\frac12)-(N+2)=s-1$). Thus $\mathcal P_0(\varphi\zeta)=0$ for $\mathcal P_0=\mathcal P-c$ and the products $\varphi\zeta$ of eigenfunctions and radial solutions.

\smallskip
\emph{(b) A subsolution.} Let $\vartheta(r)=\int_R^r\frac4t\bar c(t)\dd t$ and $\Psi=\varepsilon\varphi'\zeta_{\mu'}(1-\vartheta)$ on $\Omega'\cap\{r>R\}$. With $g=1-\vartheta$, $g'=-\frac4r\bar c$, a direct computation gives
\[
\mathcal P\Psi\le\varepsilon\varphi'\zeta_{\mu'}\Big[-g''-\Big(2\frac{\zeta_{\mu'}'}{\zeta_{\mu'}}+\frac{N+1}r+\frac r2\Big)g'+\bar c\,g\Big].
\]
Since $\zeta_{\mu'}'/\zeta_{\mu'}=-\frac r2+O(\frac1r)$, the $g'$-term equals $-2\bar c(1+O(r^{-2}))$; and $-g''=\frac4r\bar c'-\frac4{r^2}\bar c\le0$ for large $r$. As $0<g\le1$, the bracket is $\le-\bar c(1-O(r^{-2}))<0$ for $r\ge R$, $R$ large. So $\mathcal P\Psi\le0$.

\smallskip
\emph{(c) Comparison.} Choose $R\ge r_1$ as above with $\vartheta(\infty)<\frac12$, and $\varepsilon>0$ so small that $\Psi\le W$ on $\overline{D'}\times\{R\}$; this is possible because $W>0$ on this compact set. Let $Z=\varphi''(\zeta_{\mu''}+\varpi_{\mu''})$, positive on $\overline{\Omega'}\cap\{r\ge R\}$ after enlarging $R$; then $\mathcal PZ=cZ\ge0$. The function $v=(\Psi-W)/Z$ satisfies
\[
-\Delta v-\Big(\tfrac12\xi+\frac{2\nabla Z}{Z}\Big)\cdot\nabla v+\frac{\mathcal PZ}{Z}\,v\le0,
\]
and the zero-order coefficient $\mathcal PZ/Z=c$ is nonnegative. Hence the maximum principle applies to $v$: it has no positive interior maximum. On $\partial\Omega'\cap\{r>R\}$, $\Psi=0\le W$; on $\{r=R\}$, $v\le0$; and $v\to0$ as $r\to\infty$, because $\Psi,W=O(r^{s-1}e^{-r^2/4})$ while $Z\ge c\,r^{-2\lambda-1}$. Hence $v\le0$, i.e.\ $W\ge\Psi\ge\frac\varepsilon2\varphi'\zeta_{\mu'}$ in $\Omega'\cap\{r>R\}$.

\smallskip
\emph{(d) Conclusion.} $T$ preserves $|\xi|=r$, and $r^{1-s}e^{r^2/4}U(\xi)=g_r(\sigma_1(\xi))$ with $\sigma_1(\xi)=|y|/r$. So (c) gives, for $\omega\in D'$ and $r>R$,
\[
g_r(\sigma_1(\omega))-g_r(\sigma_1(T\omega))\ge c_0\,\varphi'(\omega)>0 ,
\]
where $\sigma_1(T\omega)>\sigma_1(\omega)$. Letting $r\to\infty$ along the sequence of Step~3 gives the same inequality for $\Phi/\sigma_1$ at a.e.\ $\omega\in D'$. So $\Phi/\sigma_1$ is not constant.

\emph{Step 5 (sign).} Let $g=\Phi/\sigma_1$, nonincreasing and not constant (case $\gamma>1$), and let $w(t)$ be the density of $\sigma_1$ on $\Sph^{N-1}_+$. $Z_3$ is a positive multiple of $t(2(\nu+2)t^2-3)$, $\nu=\frac{N-2}2$, resp.\ of $4t^3-3t$ for $N=2$. So it has a single zero $t_0\in(0,1)$, with $Z_3<0$ on $(0,t_0)$ and $Z_3>0$ on $(t_0,1)$. Since $\int tZ_3w=0$,
\[
c_3=\int_0^1t\,(g(t)-g(t_0))\,Z_3(t)\,w(t)\dd t ,
\]
and the integrand is $\le0$, and $<0$ on a set of positive measure. So $c_3<0$.
\end{proof}

\begin{theorem}[uniform convergence]\label{cor:unif}
$H(r,\cdot)\to\Phi$ uniformly on $\overline{\Sph^{N-1}_+}$ as $r\to\infty$.
\end{theorem}
\begin{proof}
Let $Q_R=\Rp\cap\{r>R\}$. For $\varphi\in C_c^\infty(Q_R)$, the identity in the proof of Lemma~\ref{lem:reg} gives $\int|\nabla\varphi|^2K\ge(\frac{R^2}{16}+\frac N4)\int\varphi^2K$. So for $R^2>16\lambda$ the form $\int|\nabla\varphi|^2K-\lambda\int\varphi^2K$ is coercive on $W^{1,2}_{0,K}(Q_R)$.

\emph{Step 1 (splitting).} The functional $\varphi\mapsto\int\eta_1^{2\alpha}h^q\varphi K$ is bounded on $W^{1,2}_{0,K}(Q_R)$: on each sphere, the Hardy argument of Lemma~\ref{lem:cmode}(iii) bounds $\int\sigma_1^{2\alpha+q}|\varphi|$ by the angular gradient, and the radial factor $r^{2\alpha+sq}e^{-qr^2/4}$ is integrable against $K$. By Lax--Milgram there is a unique $h_{\rm src}\in W^{1,2}_{0,K}(Q_R)$ with $(L_K-\lambda)h_{\rm src}=-\eta_1^{2\alpha}h^q$ weakly, and we put $h_{\rm hom}=h-h_{\rm src}$. Then $(L_K-\lambda)h_{\rm hom}=0$ in $Q_R$ with the boundary values of $h$. Write $H_{\rm src}$, $H_{\rm hom}$ for the corresponding rescalings. Testing with $h_{\rm src}^+$ and using coercivity gives $h_{\rm src}\le0$.

\emph{Step 2 (the source part is uniformly small).} Fix $\theta\in(\frac12,1)$ (then $\theta-2\le\gamma$, since $\gamma>-1$), let $a=s(s+N-2)+1$, and set $\hat w=e^{-r^2/4}r^s\,\varepsilon\,e^{-a/r^2}\sigma_1^\theta$. In the variable of $H$, $\hat W=\varepsilon e^{-a/r^2}\sigma_1^\theta$, and with $\Delta_{\Sph}\sigma_1^\theta=-\theta(1-\theta)\sigma_1^{\theta-2}(1-\sigma_1^2)-(N-1)\theta\sigma_1^\theta$,
\[
L_H\hat W\le-\frac{\hat W}{2r^2}-\frac{\varepsilon e^{-a/r^2}}{r^2}\,\theta(1-\theta)\,\sigma_1^{\theta-2}(1-\sigma_1^2)\qquad(r\ge R_a),
\]
where $L_H=\partial_r^2+(\frac{\kappa_1}r-\frac r2)\partial_r+r^{-2}(\Delta_{\Sph}+s(s+N-2))$. The source satisfies $F\le C_Fr^\kappa e^{-(q-1)r^2/4}\sigma_1^\gamma$, and $\sigma_1^\gamma\le C\sigma_1^\theta$ on $\{\sigma_1\ge\frac12\}$, $\sigma_1^\gamma\le\sigma_1^{\theta-2}$ on $\{\sigma_1<\frac12\}$. Hence the choice
\[
\varepsilon=\varepsilon(R)=C\sup_{t\ge R}t^{\kappa+2}e^{-(q-1)t^2/4}
\]
gives $L_H\hat W\le-F$, i.e.\ $(L_K-\lambda)\hat w\ge\eta_1^{2\alpha}h^q$. Since $\theta>\frac12$, $|\nabla\hat w|^2\sim\sigma_1^{2\theta-2}$ is integrable, so $\hat w\in W^{1,2}_{\rm loc}$ up to $\partial\Rp$, and the inequality $(L_K-\lambda)\hat w\ge\eta_1^{2\alpha}h^q$ holds weakly against nonnegative test functions in $W^{1,2}_{0,K}(Q_R)$ (integrate by parts on $\{\eta_1>\delta\}$ and let $\delta\to0$; the boundary term $\int_{\eta_1=\delta}\varphi\,\partial_1\hat w\,K$ tends to $0$ because $\varphi\in W^{1,2}_0$ and $\partial_1\hat w=O(\eta_1^{\theta-1})$ with $\theta-1>-\frac12$). Since $\hat w\ge0$ and $h_{\rm src}=0$ on $\partial Q_R$, $(\hat w+h_{\rm src})^-\in W^{1,2}_{0,K}(Q_R)$; testing $(L_K-\lambda)(\hat w+h_{\rm src})\ge0$ with it and using coercivity gives $h_{\rm src}\ge-\hat w$. Therefore
\[
|H_{\rm src}(r,\sigma)|\le\varepsilon(R)\,\sigma_1^\theta\le\varepsilon(R)\qquad(r>R),\qquad\varepsilon(R)\to0 .
\]

\emph{Step 3 (the homogeneous part converges uniformly).} $H_{\rm hom}=H-H_{\rm src}$ is bounded, and $H_{\rm hom}(R,\cdot)=H(R,\cdot)$ because $h_{\rm src}$ vanishes on $\{r=R\}$. So each of its modes is the bounded solution of the homogeneous equation with value $c_\ell(R)$ at $R$, namely
\[c_\ell(R)\,\tilde\upsilon_\ell(r)/\tilde\upsilon_\ell(R).\]
(Enlarge $R$ so that $\tilde\upsilon_\ell(R)>0$ for the finitely many $\ell\le s$; for $\ell>s$ this holds by Lemma~\ref{lem:cmode}.) For $\ell>s$ we have $\tilde\upsilon_\ell(r)/\tilde\upsilon_\ell(R)=V_\ell(z_r)/V_\ell(z_R)$, and

\smallskip\noindent
\emph{Claim:} if $A\ge6z$, then $V(4z)\le2^{1-B}V(z)$.

\smallskip\noindent
Indeed, on $\{X\ge A/2\}$, $\frac{1+X/(4z)}{1+X/z}\le\frac14+\frac34\cdot\frac{2z}A\le\frac12$; on $\{X<A/2\}$, $(1+\frac X{4z})^B\le(1+\frac A{8z})^B\le2^{-B}(1+\frac Az)^B\le2^{-B}V(z)$ by Jensen. Adding the two pieces proves the claim.

Since $V$ is decreasing, for $r\ge2R$ and $A_\ell\ge\frac32R^2$ we get $|c_\ell(R)|\,V_\ell(z_r)/V_\ell(z_R)\le\|H(R,\cdot)\|_{L^2}2^{1-B_\ell}$. With $\|Z_\ell\|_\infty\le C\ell^{(N-2)/2}$ and $B_\ell\ge\ell/2$, the series $\sum_\ell c_\ell(R)\frac{\tilde\upsilon_\ell(r)}{\tilde\upsilon_\ell(R)}Z_\ell$ converges absolutely and uniformly for $r\ge2R$, and each term converges as $r\to\infty$. The sum of this series is continuous and equals $H_{\rm hom}(r,\cdot)$ in $L^2$; since $H_{\rm hom}(r,\cdot)$ is continuous, they coincide. Hence $H_{\rm hom}(r,\cdot)$ converges uniformly on $\overline{\Sph^{N-1}_+}$.

\emph{Step 4.} By Steps 2--3, $\limsup_{r,r'\to\infty}\|H(r,\cdot)-H(r',\cdot)\|_\infty\le2\varepsilon(R)$ for every large $R$. So $H(r,\cdot)$ is uniformly Cauchy, and its uniform limit is the $L^2$ limit $\Phi$.
\end{proof}

\begin{remark}[what remains open]
\begin{enumerate}[label=\arabic*.,leftmargin=*]
\item \emph{Higher signs away from $q_c$.} Theorem~\ref{thm:sign} fixes $\operatorname{sign}c_3$ for every $q<q_c$, and Theorem~\ref{thm:nearcrit} fixes every $\operatorname{sign}c_\ell$ near $q_c$. We do not know whether, for $\gamma_c<3$, the law $\operatorname{sign}c_\ell=(-1)^{\frac{\ell-1}2}\operatorname{sign}(2\alpha+q-1)$ holds for every $q<q_c$ and every odd $\ell\ge5$.
\item Nothing is claimed on initial traces of general solutions. For $\alpha=0$, $u_s$ dominates every solution continuous on $\overline E\setminus\{(0,0)\}$ \cite[Lemma~2.7]{MV01}; that proof uses scale invariance and the Cauchy--Dirichlet problem with bounded initial data, and its extension to the weight $x_1^{2\alpha}$ is not attempted here (Theorem~\ref{thm:signed}(ii) covers self-similar solutions).
\item Sign-changing self-similar solutions for $q<q_2$: for $\alpha=0$, \cite[Remark after Theorem~2.1]{MV01} indicates that they exist for $q<(N+4)/(N+2)$ (bifurcation from higher eigenvalues of $L_K$); we have not studied this.
\end{enumerate}
\end{remark}

%---------------------------------------------------------------------
\section{Signed self-similar solutions}\label{sec:signed}

In this section $w$ may change sign. Let $\Sigma$ be the set of $w\in\X$ such that
\begin{equation}\label{eq:weaks}
\int_{\Rp}\nabla w\cdot\nabla\varphi\,K\dd\eta+\int_{\Rp}\eta_1^{2\alpha}|w|^{q-1}w\,\varphi\,K\dd\eta
=\lambda\int_Aw\varphi\,K\dd\eta\qquad\text{for all }\varphi\in\X ,
\end{equation}
and let $\mathcal S^\pm$ be the set of $w\in C^2(\Rp)\cap C(\overline{\Rp})$ that solve $L_Kw-\lambda w+\eta_1^{2\alpha}|w|^{q-1}w=0$ classically in $\Rp$, vanish on $\partial\Rp$, and satisfy $|\eta|^{2\lambda}w(\eta)\to0$ as $|\eta|\to\infty$. Then $-\Sigma=\Sigma$, $-\mathcal S^\pm=\mathcal S^\pm$, and $\Sigma\cap\{w\ge0\}=\Hh$. As in Section~\ref{sec:cont}, the profile $w=u(\cdot,1)$ of a self-similar solution $u\in C(\overline E\setminus\{(0,0)\})$ of \eqref{eq:pde}--\eqref{eq:bc} belongs to $\mathcal S^\pm$.

We use two elementary facts. For $w\in\X$, $w^\pm=\max(\pm w,0)\in\X$, with $\nabla w^+=\nabla w\,\mathbf 1_{\{w>0\}}$ a.e. And if $\varphi$ vanishes a.e.\ outside $\{w>\varepsilon\}$ for some $\varepsilon>0$, then $\nabla\varphi=0$ a.e.\ on $\{w\le\varepsilon\}$, so that, with $\mathcal B$ as in Section~\ref{sec:prelim},
\begin{equation}\label{eq:Bplus}
\mathcal B(w^+;\varphi)=\int\nabla w\cdot\nabla\varphi\,K+\int\eta_1^{2\alpha}|w|^{q-1}w\,\varphi K-\lambda\int w\varphi K .
\end{equation}

\begin{lemma}\label{lem:signedbdd}
Let $w\in\Sigma$. Then $w$ is bounded and $w(\eta)\to0$ as $|\eta|\to\infty$.
\end{lemma}
\begin{proof}
Apply the proof of Lemma~\ref{lem:reg} to $v^+$, $v=K^{1/2}w$, with the test functions $\tilde v=\zeta\,T_\tau(v^+)$. They vanish where $w\le0$, so the absorption term $\int\eta_1^{2\alpha}K^{\frac{1-q}2}|v|^{q-1}v\,\tilde v$ is nonnegative and can be dropped as before. This gives $-\Delta\widehat{v^+}\le\lambda\widehat{v^+}$ in $\mathcal D'(\R^N)$ and the bound of Lemma~\ref{lem:reg} for $w^+$. The same applies to $w^-$, since $-w\in\Sigma$.
\end{proof}

\begin{proposition}\label{prop:Spm}
$\mathcal S^\pm\subset\Sigma$ for every $q>1$.
\end{proposition}
\begin{proof}
Let $w\in\mathcal S^\pm$. We repeat the proof of Proposition~\ref{prop:S}, applying the comparison steps to $w^+$ and to $w^-=(-w)^+$.

\emph{Step 1: $|w|\le C_R\eta_1$ on $\{|\eta|\le R\}$.} In the proof of Lemma~\ref{lem:bdry} take $m=\sup\{|w(\eta)|:\eta\in\Rp,\ |\eta|\le R+2\}$ and apply Lemma~\ref{lem:comp} in $Q$ with $u_1=w^+$, $u_2=\Theta$. Each set $S_\varepsilon=\{w^+>\Theta+\varepsilon\}$ is relatively compact in $Q$ and its closure lies in the open set $G=\{w>\varepsilon/2\}$, on which $w^+=w$ is a classical solution of $\mathcal M(w)=0$. So for $\varphi$ vanishing outside $S_\varepsilon$, Lemma~\ref{lem:test}(a) (applied in $G$) gives $\mathcal B(w^+;\varphi)=0$, which is (H1). The rest of the proof is unchanged, and gives $w^+\le C_R\eta_1$. The same holds for $w^-$.

\emph{Step 2: $|w|\le CU$ in $\{r>R_0\}$.} The comparison of Step~2 of Proposition~\ref{prop:S} with $u_1=w^\pm$ and $u_2=V_\varepsilon$ works in the same way, because the sets $S_{\varepsilon'}$ are relatively compact in $Q$ and contained in $\{w^\pm>\varepsilon'\}$.

\emph{Step 3: energy and the identity \eqref{eq:weaks}.} Step~3 of Proposition~\ref{prop:S} uses the test function $w\chi_k^2\theta_j(\eta_1)^2$ and the bounds $|w|\le C_{k+1}\eta_1$, $|w|\le CU$; the sign of $w$ plays no role, and the proof of Lemma~\ref{lem:test}(a) does not use the sign of $h$ either. It gives $\nabla w\in L^2_K$ and $\eta_1^{2\alpha}|w|^{q+1}\in L^1(K)$. The functions $\big((w-\varepsilon)_+-(w+\varepsilon)_-\big)\chi_k$ have compact support in $\Rp$ and converge to $w$ in $W^{1,2}_K$, so $w\in W^{1,2}_{0,K}(\Rp)$ and $w\in\X$. Finally the proof of Lemma~\ref{lem:test}(c) is also sign-free and gives \eqref{eq:weaks} for every $\varphi\in\X$.
\end{proof}

\begin{proposition}[nonexistence]\label{prop:nonexs}
If $q\ge q_c$, then $\Sigma=\{0\}$.
\end{proposition}
\begin{proof}
Take $\varphi=w$ in \eqref{eq:weaks} and use Lemma~\ref{lem:spec}(i):
$\lambda_1\int w^2K+\int\eta_1^{2\alpha}|w|^{q+1}K\le\int|\nabla w|^2K+\int\eta_1^{2\alpha}|w|^{q+1}K=\lambda\int w^2K\le\lambda_1\int w^2K$.
Hence $\int\eta_1^{2\alpha}|w|^{q+1}K=0$ and $w=0$.
\end{proof}

\begin{proposition}[domination]\label{prop:dom}
If $1<q<q_c$, then $|w|\le h$ in $\Rp$ for every $w\in\Sigma$.
\end{proposition}
\begin{proof}
Apply Lemma~\ref{lem:comp} in $Q=\Rp$ with $u_1=w^+$ and $u_2=h$. By Lemma~\ref{lem:signedbdd}, $u_1$ is bounded and tends to $0$ at infinity, so each $S_\varepsilon$ is bounded; $u_2=h>0$. (H2) holds because $w\in\X$. In (H3), $(w^+-h-\varepsilon)_+\le(w^+-\varepsilon)_+\in W^{1,2}_0(\Rp)$, as in Proposition~\ref{prop:uniq}. For (H1): the test functions $w_1,w_2$ of the proof of Lemma~\ref{lem:comp} are bounded by a multiple of $w^+$, lie in $W^{1,2}_0(\Rp)$ and have bounded support, so they belong to $\X$; they vanish outside $\{w>\varepsilon/2\}$, so by \eqref{eq:Bplus} and \eqref{eq:weaks}, $\mathcal B(w^+;w_i)=0$, and $\mathcal B(h;w_i)=0$ by \eqref{eq:weak}. Hence $w^+\le h$. Since $-w\in\Sigma$, also $w^-\le h$.
\end{proof}

\begin{proposition}[constant sign]\label{prop:sign}
If $w\in\Sigma$ and both $w^+$ and $w^-$ are nonzero, then $\lambda>\lambda_2$, i.e.\ $q<q_2$.
\end{proposition}
\begin{proof}
Put $p=\eta_1^{2\alpha}|w|^{q-1}\ge0$ and $\|f\|^2=\int f^2K$. Testing \eqref{eq:weaks} with $\varphi=w^+$ and with $\varphi=w^-$ (both in $\X$) gives, since $\nabla w\cdot\nabla w^\pm=\pm|\nabla w^\pm|^2$ and $|w|^{q-1}w\,w^\pm=\pm p\,(w^\pm)^2$,
\[
\int|\nabla w^\pm|^2K+\int p\,(w^\pm)^2K=\lambda\|w^\pm\|^2 ,
\]
where all integrals are finite ($\int p(w^\pm)^2K=\int\eta_1^{2\alpha}(w^\pm)^{q+1}K$). The functions $w^+$ and $w^-$ have disjoint supports and $\nabla w^+\cdot\nabla w^-=0$ a.e., so for $v=aw^++bw^-$,
\[
\int|\nabla v|^2K+\int p\,v^2K=\lambda\|v\|^2 .
\]
Choose $(a,b)\ne(0,0)$ with $\int v\psi_1K=0$. Lemma~\ref{lem:spec}(ii) gives $\lambda_2\|v\|^2\le\int|\nabla v|^2K=\lambda\|v\|^2-\int p\,v^2K$. If $\lambda\le\lambda_2$, then $\int p\,v^2K\le0$, so $p\,v^2=0$ a.e. But $p>0$ on $\{w\ne0\}$, and $v=aw^+$ on $\{w>0\}$, $v=bw^-$ on $\{w<0\}$, both sets of positive measure; hence $a=b=0$, a contradiction. Therefore $\lambda>\lambda_2$.
\end{proof}

\begin{proof}[Proof of Theorem~\ref{thm:signed}]
(i) is Proposition~\ref{prop:Spm} (with the remark before Lemma~\ref{lem:signedbdd}) and Proposition~\ref{prop:nonexs}. (ii) is Proposition~\ref{prop:dom}; for $u=t^{-\lambda}w(x/\sqrt t)$ it gives $|u|\le u_s$. (iii) The first claim is Proposition~\ref{prop:sign}. If $q_2\le q<q_c$ and $w\in\Sigma\setminus\{0\}$, then $w$ or $-w$ is nonnegative, hence in $\Hh\setminus\{0\}$, hence equal to $h$ by Proposition~\ref{prop:uniq}. The statement on continuous solutions follows from (i).
\end{proof}

\begin{remark}\label{rem:signed}
For $\alpha=0$, Proposition~\ref{prop:dom} is \cite[Lemma~2.5, Step~5]{MV01}, and Theorem~\ref{thm:signed}(iii) is \cite[Lemma~2.6]{MV01} for $(N+4)/(N+2)<q<q_c$. There the operator $L_K+|w|^{q-1}$ is used, with a bounded potential. The argument above uses only $w^\pm$, so it includes the endpoint $q=q_2$ and allows the singular weights $\alpha<0$, for which $\eta_1^{2\alpha}|w|^{q-1}$ may be unbounded near $\partial\Rp$.
\end{remark}

\section*{Statement on the use of AI tools}
AI-assisted tools were used exclusively for language editing and grammar correction.
%---------------------------------------------------------------------

\end{document}